\documentclass[12pt,draft]{extartic}

\usepackage{amsmath,amsthm,amssymb,calc}
\usepackage[dvips]{graphics, graphicx}

\usepackage[english]{babel}                         

\newcounter{num}[section]

\newcommand{\Th}{\refstepcounter{num}
{\bf Theorem \arabic{section}.\arabic{num} }}

\newcommand{\Lemma}{\refstepcounter{num}
{\bf Lemma \arabic{section}.\arabic{num} }}

\newcommand{\Pred}{\refstepcounter{num}
{\bf Proposition \arabic{section}.\arabic{num} }}

\newcommand{\Note}{\refstepcounter{num}
{\it Note \arabic{section}.\arabic{num} }}

\newcommand{\Exm}{\refstepcounter{num}
{\bf Example \arabic{section}.\arabic{num} }}

\newcommand{\St}{\refstepcounter{num}
{\bf Statement \arabic{section}.\arabic{num} }}

\newcommand{\Def}{\refstepcounter{num}
{\it Definition \arabic{section}.\arabic{num} }}

\newcommand{\Proof}{{\bf Proof. }}

\def\eps{\varepsilon}
\def\_phi{\varphi}

\def\a{\alpha}
\def\d{\delta}
\def\l{\lambda}

\def\v{\vec}
\def\F{\widehat}
\def\L{\Lambda}
\def\m{\times}
\def\t{\tilde}
\def\o{\omega}

\def\z{{\mathbb Z}}
\def\C{{\mathbb C}}
\def\r{\mathcal{R}}
\def\Z_N{{\mathbb Z}_N}
\def\Z{{\mathbb Z}}

\author{Shkredov I.D.}
\title{Some examples of sets of large exponential sums
\footnote{This work was supported by RFFI grant no. 06-01-00383,
President's of Russian Federation grant N 1726.2006.1 and INTAS
(grant no. 03--51--5-70).}}
\date{}
\begin{document}
\maketitle

\begin{center}
    Annotation.
\end{center}

{\it \small Let $A$ be a subset of $\mathbb{Z} / N\mathbb{Z}$ and
let $\mathcal{R}$ be the set of large Fourier coefficients of $A$.
Properties of $\mathcal{R}$ have been  studied in works of
M.--\,C. Chang, B. Green and the author. In the paper we obtain
some new results on sets of large exponential sums.
}
\\
\\
\\

\refstepcounter{section}

{\bf \arabic{section}. Introduction.}

Let $N$ be a positive integer. By $\Z_N$ denote the set $\z / N
\z$. Let $f: \Z_N \to \C$ be an arbitrary function. Denote by
$\F{f}$ the Fourier transform of $f$
\begin{equation}\label{}
    \F{f} (r) = \sum_{n\in \Z_N} f(n) e(-nr) \,,
\end{equation}
where $e(x) = e^{-2\pi i x/N}$.

 Let $\d,\a$ be real numbers, $0<\a \le \d \le 1$ and
 let $A$ be a subset of $\Z_N$ of cardinality $\d N$.
 It is very convenient to write $A(x)$ for such a function.
 Thus $A(x) = 1$ if $x\in A$ and $A(x) = 0$ otherwise.
 Consider the set $\r_\a$ of large exponential sums of the set $A$
\begin{equation}\label{f:R_def}
    \r_\a = \r_\a (A) = \{~ r\in \Z_N ~:~ |\F{A} (r)| \ge \a N ~\} \,.
\end{equation}
 In many problems of combinatorial number theory is important to know the structure of the set $\r_\a$
 (see \cite{Gow_surv}).
 In other words what kind of properties $\r_\a$ has?

 In 2002 M.--\,C. Chang proved the following result \cite{Ch_Fr}.

 \Th {\bf (Chang)}
 \label{t:Chang}
 { \it
    Let $\d,\a$ be real numbers, $0<\a \le \d \le 1$,
    $A$ be a subset of $\Z_N$, $|A| = \d N$.
    Then there exists a set $\L = \{ \lambda_1, \dots, \lambda_{|\L|} \} \subseteq \Z_N$,
    $|\L| \le 2 (\d / \a)^2 \log ( 1/\d)$ such that
    for any $r \in \r_\a$ we have
    \begin{equation}\label{f:presentation}
        r = \sum_{i=1}^{|\L|} \eps_i \lambda_i  \pmod N \,,
    \end{equation}
    where
    $\eps_i \in \{ -1,0,1 \}$.
 }

 Using approach of paper \cite{Ruzsa_Freiman} (see also  \cite{Bilu})
 Chang applied her result to prove the famous Freiman's theorem \cite{Freiman}
 on sets with small doubling.
 Another applications of Theorem \ref{t:Chang} were obtained by B. Green in \cite{GreenA+A},
 and by T. Schoen in \cite{Schoen}.
%
 If the parameter $\a$ is close to $\d$ then the structural properties of the set $\r_\a$
 was studied in papers
 \cite{Freiman_Yudin,Besser,Lev}, see also survey \cite{Konyagin_Lev}.


 In paper \cite{Green_Chang_exact} Green showed that Chang's theorem is sharp in a certain sense.
 Let $E = \{ e_1, \dots, e_{|E|} \} \subseteq \Z_N$ be an arbitrary set.
 By ${\rm Span}(E)$ denote the set of all sums $\sum_{i=1}^{|E|} \eps_i e_i$, where
 $\eps_i \in \{ -1,0,1 \}$.

 \Th {\bf (Green)}
 \label{t:Green_Chang_exact}
 { \it
    Let $\d,\a$ be real numbers, $\d \le 1/8$, $0< \a \le \d/32$.
    Suppose that
    \begin{equation}\label{}
        \left( \frac{\d}{\a} \right)^2 \log (1/ \d) \le \frac{\log N}{\log \log N} \,.
    \end{equation}
    Then there exists a set $A\subseteq \Z_N$, $|A| = [ \d N ]$ such that
    the set $\r_\a$ does not contain in ${\rm Span}(\L)$
    for any set
    $\L$ of cardinality $2^{-12} (\d/ \a)^2 \log (1/ \d)$.
 }

 In papers \cite{Sh_dokl_exp1,Sh_exp1} further results on sets of large
 exponential sums were obtained.
 In particular the author proved the following theorem

 \Th
 \label{t:main}
 {\it
    Let $\d,\a$ be real numbers, $0< \a \le \d$,
    $A$ be a subset of $\Z_N$, $|A| = \d N$, and
    $k\ge 2$ be a positive integer.
    Let also $B\subseteq \r_\a \setminus \{ 0 \}$ be an arbitrary set.
    Then the number
    \begin{equation}\label{f:T_k_def}
        T_k (B) := |\{ ~ (r_1,\dots, r_k, r_1', \dots, r_k') \in B^{2k} ~:~
                    r_1 + \dots + r_k = r_1' + \dots + r_k' ~ \}|
    \end{equation}
    is at least
    \begin{equation}\label{f:T_k}
        \frac{\d \a^{2k}}{2^{4k} \d^{2k}} |B|^{2k} \,.
    \end{equation}
 }

 In article \cite{Sh_exp1} was showed that
 Theorem \ref{t:main} and an inequality of W. Rudin \cite{Rudin_book}
 imply M.--\,C. Chang's theorem.
 Moreover in \cite{Sh_exp1} the following
 improvement of Theorem \ref{t:Chang} was obtained.

 \Th
   \label{t:Chang_log}
   {\it
    Let $N$ be a positive integer, $(N,6)=1$,
    $\d,\a$ be real numbers, $0< \a \le \d \le 1/16$,
    and
    $A$ be a subset of $\Z_N$, $|A| = \d N$.
    Then there exists a set
    $\L^* \subseteq \Z_N$,
    \begin{equation}\label{f:L_est}
        |\L^*| \le \min \,( \, \max(\, 2^{30} (\d /\a)^2 \log (1/\d), 2^{ 4 ( \log \log (1/\d) )^2 + 2} \, ),
                                            \,\,  2^{20} (\d/\a)^2 \log^{13/7} (1/\d) \,)
    \end{equation}
    such that
    for any
    $r \in \r_\a$ there exists a tuple $\l_1^*, \dots, \l_M^* \in \L^*$,
    $M \le s8 \log (1/\d)$ such that
    \begin{equation}\label{f:r=log}
        r = \sum_{i=1}^M \eps_i \l^*_i \pmod N \,,
    \end{equation}
    where
    $\eps_i \in \{ -1,0,1 \}$.

     Besides there exists a set $\t{\L} \subseteq \Z_N$,
    \begin{equation}\label{f:tL_est}
        |\t{\L}| \le 2^{20} (\d/\a)^2 \log^{5/3} (1/\d) \log \log (1/\d)
    \end{equation}
    such that for any residual
    $r \in \r_\a$ there exists a tuple $\t{\l}_1, \dots, \t{\l}_M \in \t{\L}$, $M\le 8 \log (1/\d)$
     such that (\ref{f:r=log}) holds.
   }

    The paper is organized as follows.

    In section \ref{examples}
    we show that Theorem \ref{t:main}
    is sharp in a certain sense.
    In our proof we construct concrete sets $A \subseteq \Z_N$
    with required properties.
    Besides in the section we obtain a result which is an inverse to Chang's theorem in some sense
    (see Theorem \ref{t:Green+}).

    In \S  \ref{positive} we obtain the following improvement of Theorem \ref{t:Chang_log}
    (we consider the case when $N$ is a prime number).

    \Th
   \label{t:Chang_log_m}
   {\it
    Let $N,d$ be positive integers,
    $\d,\a$ be real numbers, $0< \a \le \d \le 2^{-8}$,
    $\_phi = (\sqrt{73} - 5 ) / 2$,
    and
    $A$ be a subset of $\Z_N$, $|A| = \d N$.
    Then there exists a set
    $\L^* \subseteq \Z_N$,
    \begin{equation}\label{f:L_est+}
        \min \,(
       \, \max(\, 2^{30 + 8d (\log (1/\d))^{-1}} \cdot (\d /\a)^2 \log (1/\d), 2^{( \log \log (1/\d) )^2 + 3} \, ),
                                            \,\,  2^{20 + 8d (\log (1/\d))^{-1}} (\d/\a)^2 \log^{\_phi} (1/\d) \, )
    \end{equation}
    such that
    for any residual
    $r \in \r_\a \setminus \{ 0 \}$ there exists a matrix
    $M = (m_{ij})_{i\in [d], j\in [|\L^*|]}$ of rank $d$
    such that for any $i\in [d]$ we have $\sum_{j=1}^{|\L^*|} |m_{ij}| \le 4 \log (1/\d)$
    and for all  $i\in [d]$
    \begin{equation}\label{f:r=log+}
        r = \sum_{j=1}^{|\L^*|} m_{ij} \l^*_j \pmod N \,.
    \end{equation}
   }

   Certainly, our question on the structure of $\r_\a$ (as any question of combinatorial number theory)
   can be asked for
   any finite Abelian group $G$ not only for $\Z_N$.
   It turns out that (see \cite{Green_Sz_Ab,Bu} and, especially, a wonderful survey \cite{Green_finite_fields})
   many problems of combinatorial number theory are considerably easier in groups
   $\Z_p^n$, where $p$ is a small prime number
   (for example $p=2,3$ or $5$).

   In section \ref{finite_fields} we obtain some analogs of results of section \ref{examples}
   for groups $\Z_p^n$.
   Main ideas of the proofs are easier in the groups then in $\Z_N$,
   and all technical details are simplified.

    Finally note that the statements of Theorems \ref{t:Chang}, \ref{t:Chang_log}, and \ref{t:Chang_log_m}
    are trivial if the parameter $\d$ does not tends to zero as $N \to \infty$.
    In the case there are not non--trivial
    restrictions on structure of the set $\r_\a$
    (see papers \cite{Katznelson,Nazarov,Ball}).

In our forthcoming papers we are going to obtain further results
on sets of large exponential sums.

The author is grateful to Professor N.G. Moshchevitin for constant
attention to this work.

\refstepcounter{section} \label{examples}

{\bf \arabic{section}. Some examples of sets of large exponential
sums.}

Let $N$ be a positive integer. It is very convenient to write
$[N$] for $\{1,2,\dots, N\}$. Let also $f: \Z_N \to \C$ be an
arbitrary function. By Parseval's identity
\begin{equation}\label{f:Par}
    \sum_{r\in \Z_N} |\F{f} (r)|^2  = N \sum_{n\in \Z_N} |f(n)|^2 \,.
\end{equation}

Let $\d,\a$ be real numbers, $0<\a \le \d \le 1$ and let $A$ be a
subset of $\Z_N$, $|A| = \d N$.
 It is easy to see that $0\in \r_\a$ and $\r_\a = -\r_\a$.
 Further, using (\ref{f:Par}), we obtain $|\r_\a| \le \d / \a^2$.

In the section we give some examples of sets of large exponential
sums. First of all note that any "small"\, subset of $\Z_N$ is a
set of large exponential sums. To be precise, we have the
following proposition.

\Pred \label{pr:small} {\it
    Let $\d,\a \in (0,1]$ be real numbers, $\d \le 1/2$, $20 N^{-1/2} < \a \le \d/2$,
    and $S \subseteq \Z_N$ be an arbitrary set such that $0\in S$, $S=-S$ and
    $|S| \le \d / (2\a)$.
    Then there exists a set $A\subseteq \Z_N$, $|A| = [\d N]$ such that $\r_\a (A) = S$.
}


    To prove Proposition \ref{pr:small} we need in the following well--known lemma
    (see \cite{Spencer_Beck} and \cite{Green_Chang_exact}).

\Lemma \label{l:char_f} {\it
    Let $f : \Z_N \to [0,1]$ be a function.
    Then there exists a set $C \subseteq \Z_N$, $|C| = [\sum_{x\in \Z_N} f(x)]$
    such that for all $r\in \Z_N \setminus \{ 0 \}$, we have
    $|\F{C} (r) - \F{f} (r)| \le 20 \sqrt{N}$.
}

{\bf Proof of Proposition \ref{pr:small}. }
    Let $S^* = S \setminus \{ 0 \}$.
    Consider the function $f(x) = \d + 2\a \sum_{r\in S^*} e(rx)$.
    Since $S=-S$, it follows that $f(x)$ is a real function.
    We have $|S| \le \d / (2\a)$ and $\d \le 1/2$.
    Hence for all $x\in \Z_N$, we get $0\le f(x) \le 1$.
    Besides $\sum_{x\in \Z_N} f(x) = \d N$ and for any $r\in \Z_N \setminus \{ 0 \}$, we have
    $\F{f} (r) = 2\a S^* (r) N$.
    Using Lemma \ref{l:char_f}, we obtain the set $A$ such that
    $|A| = [\sum_{x\in \Z_N} f(x)] = [\d N]$ and for all $r\in \Z_N \setminus \{ 0 \}$
    $$
        | \F{A} (r) - \F{f} (r) | = | \F{A} (r) - 2\a S^* (r) N | \le 20 \sqrt{N} \,.
    $$
    Since $\a > 20 N^{-1/2}$ it follows that for any $r\in S^*$ the following inequality holds
    $|\F{A} (r)| \ge 2\a N - 20 \sqrt{N} \ge \a N$.
    Hence $S \subseteq \r_\a (A)$.
    Using the inequality $\a > 20 N^{-1/2}$ again, we obtain
    $|\F{A} (r)| < \a N$ for all $r\notin S^*$, $r\neq 0$.
    Whence $\r_\a (A) = S$.
    This completes the proof.

So any small symmetrical subset of $\Z_N$ is a set of large
exponential sums.
What is the structure of the large sets $\r_\a$?
This question is not easy but
clearly, these sets have special properties. For example any set
$\r_\a$ has large quantity $T_k (\r_\a)$ (see Theorem
\ref{t:main}).

Chang's theorem is another example which is demonstrating that our
sets have really specific properties. This theorem can be
reformulate as follows : any set of large exponential sums has
small dissociated subset. We say that a set $\mathcal{D} = \{ d_1,
\dots, d_{|\mathcal{D}|} \} \subseteq \Z_N$ is {\it dissociated}
if the equality
   \begin{equation}\label{f:def_diss}
        \sum_{i=1}^{|\mathcal{D}|} \eps_i d_i = 0 \pmod N \,,
   \end{equation}
where $\eps_i \in \{ -1,0,1 \}$ implies that all $\eps_i$ are
equal to zero. Actually Chang proved that {\it any} dissociated
subset of $\r_\a (A)$ has cardinality less than $2 (\d / \a)^2
\log ( 1/\d)$. Now if we let $\L$ be a maximal dissociated subset
of $\r_\a (A)$ then it is easy to see that for any $r\in \r_\a$ we
have (\ref{f:presentation}) (see details in \cite{Ch_Fr} or in
\cite{Green_Chang2}).

In the section we obtain a result which is an inverse to Chang's
theorem in some sense. We show that any not very big dissociated
subset of $\Z_N$ is a set of large exponential sums. We
extensively use approach of B. Green (see
\cite{Green_Chang_exact}) in our prove. His method is connected
with "niveau sets"\, of I. Ruzsa (see \cite{RuzsaA+A}).

Let us introduce a few further pieces of notation.
We say that a set $\mathcal{D} \subseteq \Z_N$ is {\it
$k$---dissociated} if the equality (\ref{f:def_diss}), where
$|\eps_i| \le k$ implies that all $\eps_i$ are equal to zero.
Using this definition we can reformulate Theorem
\ref{t:Green_Chang_exact} as follows.

\Th {\bf (Green)} \label{t:Green_actual} {\it
    Let  $\d,\a$ be real numbers, $\d \le 1/8$, $20 N^{-1/2} < \a \le \d /32$,
    and $\L$ be a $6|\L|$---dissociated set, $|\L| \le 2^{-11} (\d/\a)^2 \log (1/\d)$.
    Then there exists a set $A\subseteq \Z_N$, $|A| = [\d N]$ such that
    $\r_\a (A) = \{ 0 \} \bigsqcup \L \bigsqcup -\L$.
}

\Note
If $\L$ is a $6 |\L|$ --- dissociated set then we have $|\L| \ll
\log N / \log \log N$. Thus Theorem \ref{t:Green_actual} is not
useful for sets of cardinality $ \gg \log N / \log \log N$.

We do not prove Theorem \ref{t:Green_actual} here and obtain a
slightly stronger result. To formulate this result we need in the
following definition (see \cite{Ruzsa_independent} and
\cite{Sh_dokl_exp1,Sh_exp1}).

\Def \label{d:L_k,s}
    Let $k,s$ be positive integers.
    Consider the family $\L(k,s)$ of subsets of $\Z_N$.
    A set $\L = \{ \l_1, \dots, \l_{|\L|} \}$ belongs to the family
   $\L(k,s)$ if the equality
   \begin{equation}\label{f:d_L_k,s}
        \sum_{i=1}^{|\L|} \l_i s_i = 0 \pmod N \,, \quad \l_i \in \L \,, \quad s_i \in \Z \,, \quad
        |s_i| \le s\,, \quad \sum_{i=1}^{|\L|} |s_i| \le k \,,
   \end{equation}
   implies that $s_i$ are equal to zero.

    A set $\L = \{ \l_1, \dots, \l_{|\L|} \}$ belongs to the family
   $\L(k,\infty)$
   if the equality
   \begin{equation}\label{}
        \sum_{i=1}^{|\L|} \l_i s_i = 0 \pmod N \,, \quad \l_i \in \L \,, \quad s_i \in \Z \,, \quad
        \quad \sum_{i=1}^{|\L|} |s_i| \le k \,,
   \end{equation}
   implies that all $s_i$ are equal to zero.


    Note that for any $\L \in \L(k,s)$, where $s$ is a positive integer or $\infty$,
    we have $0 \notin \L$ and $\L \cap -\L = \emptyset$.

    We need in a more delicate definition.

\Def Let $k,p$ be positive integers and $s$ be a positive integer
or $\infty$.
   Consider the family
   $\L(k,s,p)$ of subsets of $\Z_N$.
   A set $\L = \{ \l_1, \dots, \l_{|\L|} \}$ belongs to the family
   $\L(k,s,p)$ if there exists a partition of $\L$ into $p$ subsets
   $\L_1 = \{ \l^{(1)}_1, \dots, \l^{(1)}_{|\L_1|} \}, \dots, \L_p = \{ \l^{(p)}_1, \dots, \l^{(p)}_{|\L_p|} \}$,
   such that the cardinalities of any two subsets $\L_i$, $\L_j$ differ by at most two times,
   and the equality
   \begin{equation}\label{}
        \sum_{i=1}^{|\L_1|} \l^{(1)}_i s^{(1)}_i + \dots + \sum_{i=1}^{|\L_p|} \l^{(p)}_i s^{(p)}_i= 0 \pmod N \,, \quad \mbox{ where }
   \end{equation}
   \begin{equation}\label{}
        \quad \l^{(i)}_j \in \L_i \,, \quad s^{(i)}_j \in \Z \,, \quad
        |s^{(i)}_j| \le s\,, \quad \quad \sum_{j=1}^{|\L_i|} |s^{(i)}_j| \le k \,, \quad i=1,\dots, p
   \end{equation}
   imply that all $s^{(i)}_j$, $i=1,\dots, p$,\, $j = 1,\dots, |\L_i|$ are equal to zero.

\Exm Let $k,s,p$ be positive integers. Then any set $\L \in
\L(kp,s)$, $|\L| \ge p$ belongs to the family $\L(k,s,p)$.

\Th \label{t:Green+} {\it
    Let $\d,\a$ be real numbers, $\d \le 1/8$, $640 N^{-1/2} < \a \le 2^{-27} \d$,
    and $\L$ be a subset of the family $\L( (\d/\a)^2, (\d/\a)^2, [\log (1/\d)])$,
    $|\L| \le 2^{-12} (\d/\a)^2 \log (1/\d)$.
    Then there exists a set $A\subseteq \Z_N$, $|A| = [\d N]$ such that
    $\r_\a (A) = \{ 0 \} \bigsqcup \L \bigsqcup -\L$.
}

    To prove Theorem \ref{t:Green+} we need in the following auxiliary result.

\St \label{st:small_diss} {\it
    Let $\d,\a \in (0,1]$ be real numbers, $640 N^{-1/2} < \a \le 2^{-10} \d$,
    and $\L$ be any 2---dissociated set of the cardinality at most $\frac{\d}{3\a} \log (1/\d)$.
    Then there exists a set $A\subseteq \Z_N$, $|A| = [\d N]$ such that
    $\r_\a (A) = \{ 0 \} \bigsqcup \L \bigsqcup -\L$.
}

    Thus Statement \ref{st:small_diss} tells us that any
    2---dissociated set  of the cardinality rather more then $\d/\a$                                    
    is a set of large exponential sums.

{\bf Proof of Statement \ref{st:small_diss}. }
    Let $\L = \{ \l_1, \dots, \l_{|\L|} \}$ be a dissociated set,
    $|\L| \le \frac{\d}{3\a} \log (1/\d)$.
    Let also $m=|\L|$ and $c=(3 \ln 2) /2$.
    Consider the function
    \begin{equation}\label{tmp:23:04:15}
        f(x) = \d \prod_{j=1}^{m} ( 1 + \frac{2c\a}{\d} \cos (\l_j x) )
                    = \d \prod_{j=1}^{m} ( 1 + \frac{c\a}{\d} ( e(\l_j x) + e(-\l_j x)) ) \,.
    \end{equation}
    Clearly, $f(x) \ge 0$ and $\sum_{x\in \Z_N} f(x) = \d N$.
    By assumption  $m \le \frac{\d}{3\a} \log (1/\d)$.
    Hence $f(x) \le \d (1+\frac{2c\a}{\d})^m \le 1$.
    Let
    $$
        \nu_m (n) = | \{ r_1, \dots, r_m \in \L ~:~ n = \pm r_1 \pm \dots \pm r_m \} |\,, \dots ,
    $$
    $$
        \nu_2 (n) = | \{ r_1, r_2 \in \L ~:~ n = \pm r_1 \pm r_2 \} |\,, \quad
        \nu_1 (n) = | \{ r_1 \in \L ~:~ n = \pm r_1 \} | \,.
    $$
    So the quantity $\nu_s (n)$ is the number of the solutions of the equation
    $n = \sum_{i=1}^s \eps_i r_i$, where $\eps_i = \pm 1$, and $r_i \in \L$.
    Using (\ref{tmp:23:04:15}), we get
    \begin{equation}\label{}
        f(x) = \d + \d \frac{c\a}{\d} \sum_n \nu_1 (n) e(nx) + \d \left( \frac{c\a}{\d} \right)^2 \sum_n \nu_2 (n) e(nx) +
        \dots + \d \left( \frac{c\a}{\d} \right)^m \sum_n \nu_m (n) e(nx) \,.
    \end{equation}
    In other words
    \begin{equation}\label{}
        f(x) = \d + c\a \F{\nu}_1 (-x) + \d \left( \frac{c\a}{\d} \right)^2 \F{\nu}_2 (-x) +
        \dots + \d \left( \frac{c\a}{\d} \right)^m \F{\nu}_m (-x) \,.
    \end{equation}
    Hence
    \begin{displaymath}
        \F{f} (r) = \left\{ \begin{array}{ll}
                        \d N \,,   & \mbox{ if } r = 0, \\
                        N ( c\a \cdot \nu_1 (r) + \d \left( \frac{c\a}{\d} \right)^2 \nu_2 (r) + \dots
                                + \d \left( \frac{c\a}{\d} \right)^m \nu_m (r) )  \,,   & \mbox{ otherwise. }
                \end{array} \right.
    \end{displaymath}
    By assumption $\L$ is a 2---dissociated set.
    It follows that for all $i \ge 1$, we have $\nu_i (n) \le 1$.
    If $r \in \L$ or $r\in -\L$ then
    \begin{equation}\label{TMP1}
        |\F{f} (r)| \ge c\a N -
                                N \left( \d \left( \frac{c\a}{\d} \right)^2 + \d \left( \frac{c\a}{\d} \right)^3 + \dots
                                    + \d \left( \frac{c\a}{\d} \right)^m \right) \ge (1+2^{-5}) \a N \,.
    \end{equation}
    Similarly if $r \notin \L \bigsqcup -\L$ then
    \begin{equation}\label{TMP2}
        |\F{f} (r)| \le N \left( \d \left( \frac{c\a}{\d} \right)^2 + \d \left( \frac{c\a}{\d} \right)^3 + \dots
                                    + \d \left( \frac{c\a}{\d} \right)^m \right) \le \frac{1}{2} \a N \,.
    \end{equation}
    Using Lemma \ref{l:char_f}, we find a set $A$ such that
    $|A| = [\sum_{x\in \Z_N} f(x)] = [\d N]$ and for all $r\in \Z_N \setminus \{ 0 \}$, we have
    $
        | \F{A} (r) - \F{f} (r) | \le 20 \sqrt{N} \,.
    $
    Combining (\ref{TMP1}), (\ref{TMP2}) and $\a > 640 N^{-1/2}$,
    we get
    $\r_\a (A) = \{ 0 \} \bigsqcup \L \bigsqcup -\L$.
    This completes the proof of Statement \ref{st:small_diss}.

    Let us return to the proof of Theorem \ref{t:Green+}.

    We need in a lemma from \cite{Green_Chang_exact}.

\Lemma \label{l:p_k} {\it
    Let $k$ be a positive integer, and
    \begin{equation}\label{}
        p_k (x) = 2 + x \sum_{j=0}^k \frac{(-1)^j x^{2j}}{2^{4j} j!} \,.
    \end{equation}
    Then for all $x$ such that $|x| \le \sqrt{k}$, we have $0\le p_k(x) \le 4$.
}

{\bf Proof of Theorem \ref{t:Green+}}.
    Let $k = s = (\d/\a)^2$, $p = [\log (1/\d)]$.
    Let also $k_i = |\L_i|$ and $\L_i = \{ \l^{(i)}_1, \dots, \l^{(i)}_{k_i} \}$, $i=1,2,\dots, p$.
    By definition of the family $\L(k,k,p)$, we get
\begin{equation}\label{tmp:five}
    \frac{|\L|}{2p} \le k_i \le \frac{2|\L|}{p} \,.
\end{equation}
    We can assume that
\begin{equation}\label{tmp:low}
    |\L| > \frac{\d}{3\a} \log (1/\d) > 2 \log (1/\d) \,.
\end{equation}
    Indeed if $|\L| \le \d/ (3\a) \cdot \log (1/\d)$ then for all $i\in [p]$, we have
    $2 k_i + 1 \le 8\d / \a \le k$.
    Besides $s\ge 2$.
    Hence $\L$ is a 2---dissociated set,
    and using Statement \ref{st:small_diss},
    we obtain the required set $A$.

    Let
\begin{equation}\label{tmp:staR}
    g(x) = 4^{-p} \prod_{i=1}^p p_{k_i} \left( \frac{\cos(2\pi \l^{(i)}_1 x / N) + \dots + \cos(2\pi \l^{(i)}_{k_i} x / N)}
                                                    {\sqrt{k_i}} \right) \,.
\end{equation}
    Using Lemma \ref{l:p_k}, we get for all $x \in \Z_N$ the following inequality holds $0 \le g(x) \le 1$.
    Consider the $i$--th term of product (\ref{tmp:staR}).
    By formula $\cos (2\pi x/N) = (e(x) + e(-x)) /2$, we obtain
$$
    p_{k_i} \left( \frac{\cos(2\pi \l^{(i)}_1 x / N) + \dots + \cos(2\pi \l^{(i)}_{k_i} x / N)}{\sqrt{k_i}} \right)
        =
$$
\begin{equation}\label{tmp:staRR}
            2 + \frac{1}{2 \sqrt{k_i}} \sum_{j=0}^{k_i} \frac{(-1)^j}{(64 k_i)^j j!}
                \left( e(\l^{(i)}_1 x) + e(-\l^{(i)}_1 x) + \dots + e(\l^{(i)}_{k_i} x) + e(-\l^{(i)}_{k_i} x) \right)^{2j+1} \,.
\end{equation}
    Thus
$$
    g(x) = \sum_{\sum_{l=1}^{k_1} |s^{(1)}_l| \le 2k_1 +1} \dots \sum_{\sum_{l=1}^{k_p} |s^{(p)}_l| \le 2k_p +1}
$$
\begin{equation}\label{tmp:1_dashed}
                Q(s^{(1)}_1, \dots, s^{(1)}_{k_1}, \dots, s^{(p)}_1, \dots, s^{(p)}_{k_p})
                e((s^{(1)}_1 \l^{(1)}_1 + \dots + s^{(1)}_{k_1} \l^{(1)}_{k_1} +
                    \dots + s^{(p)}_1 \l^{(p)}_1 + \dots + s^{(p)}_{k_p} \l^{(p)}_{k_p}) x) \,,
\end{equation}
where $Q(s^{(1)}_1, \dots, s^{(1)}_{k_1}, \dots, s^{(p)}_1, \dots,
s^{(p)}_{k_p})$ is a coefficient which attached to $
    e((s^{(1)}_1 \l^{(1)}_1 + \dots + s^{(1)}_{k_1} \l^{(1)}_{k_1} +
                    \dots + s^{(p)}_1 \l^{(p)}_1 + \dots + s^{(p)}_{k_p} \l^{(p)}_{k_p}) x)
$.
    Similar
$$
    p_{k_i} \left( \frac{\cos(2\pi \l^{(i)}_1 x / N) + \dots + \cos(2\pi \l^{(i)}_{k_i} x / N)}{\sqrt{k_i}} \right)
        =
$$
\begin{equation}\label{}
        =
        \sum_{\sum_{l=1}^{k_i} |s^{(i)}_l| \le 2k_i +1} Q(s^{(i)}_1, \dots, s^{(i)}_{k_i})
                e((s^{(i)}_1 \l^{(i)}_1 + \dots + s^{(i)}_{k_i} \l^{(i)}_{k_i}) x)
\end{equation}
    Clearly,
    $Q(s^{(1)}_1, \dots, s^{(1)}_{k_1}, \dots, s^{(p)}_1, \dots, s^{(p)}_{k_p})
        =   4^{-p}
            \prod_{i=1}^p Q(s^{(i)}_1, \dots, s^{(i)}_{k_i})$.
    By assumption $\L$ belongs to $\L(k,k,p)$.
    This implies that all the sums
$
    s^{(1)}_1 \l^{(1)}_1 + \dots + s^{(1)}_{k_1} \l^{(1)}_{k_1} +
                    \dots + s^{(p)}_1 \l^{(p)}_1 + \dots + s^{(p)}_{k_p} \l^{(p)}_{k_p}
$ are distinct. In particular
\begin{equation}\label{tmp:(1)}
    \sum_{x\in \Z_N} g(x) = 2^{-p} N \,.
\end{equation}
Using formula (\ref{tmp:1_dashed}), we obtain that any non--zero
Fourier coefficient of the function $g(x)$ must have the following
form $
    s^{(1)}_1 \l^{(1)}_1 + \dots + s^{(1)}_{k_1} \l^{(1)}_{k_1} +
                    \dots + s^{(p)}_1 \l^{(p)}_1 + \dots + s^{(p)}_{k_p} \l^{(p)}_{k_p}
$, where $\sum_{l=1}^{k_i} |s^{(i)}_l| \le 2k_i +1$, $i\in [p]$.
Moreover, any Fourier coefficient of  $g(x)$ of the form $
    s^{(1)}_1 \l^{(1)}_1 + \dots + s^{(1)}_{k_1} \l^{(1)}_{k_1} +
                    \dots + s^{(p)}_1 \l^{(p)}_1 + \dots + s^{(p)}_{k_p} \l^{(p)}_{k_p}
$ is equal to $N \cdot Q(s^{(1)}_1, \dots, s^{(1)}_{k_1}, \dots,
s^{(p)}_1, \dots, s^{(p)}_{k_p})$. Let us prove that for all $i\in
[p]$, $j\in [k_i]$, we have
\begin{equation}\label{tmp:(2)}
    |\F{g} (\l^{(i)}_j) | = |\F{g} (-\l^{(i)}_j) | \ge 2^{-p} \frac{N}{8 \sqrt{k_i}} \,.
\end{equation}
Clearly, it suffices to deal with the case $i=j=1$. In other words
we need to find the coefficient $Q(1,0,\dots, 0)$.
We have
$$
    p_{k_1} \left( \frac{\cos(2\pi \l^{(1)}_1 x / N) + \dots + \cos(2\pi \l^{(1)}_{k_1} x / N)}{\sqrt{k_1}} \right)
        =
$$
$$
    =
          2 + \frac{1}{2 \sqrt{k_1}} \sum_{j=0}^{k_1} \frac{(-1)^j}{(64 k_1)^j j!}
                \left( e(\l^{(1)}_1 x) + e(-\l^{(1)}_1 x) + \dots + e(\l^{(1)}_{k_1} x) + e(-\l^{(1)}_{k_1} x) \right)^{2j+1}
    =
$$
\begin{equation}\label{tmp:(3)}
    =
        \sum_{\sum_{l=1}^{k_1} |s^{(1)}_l| \le 2k_1 +1} Q(s^{(1)}_1, \dots, s^{(1)}_{k_1})
                e((s^{(1)}_1 \l^{(1)}_1 + \dots + s^{(1)}_{k_1} \l^{(1)}_{k_1}) x) \,.
\end{equation}
    The coefficient attached to $e(\l^{(1)}_1)$ in (\ref{tmp:(3)}) from $j=0$ equals $1/(2\sqrt{k_1})$.
    Let us prove that the sum of the coefficients attached to $e(\l^{(1)}_1)$ from $j\ge 1$
    at most $1/(4 \sqrt{k_1})$.

    Let $l=1,2,\dots, k_1$, and consider the product of $(2l+1)$ brackets
$
    \left( e(\l^{(1)}_1 x) + e(-\l^{(1)}_1 x) + \dots + e(\l^{(1)}_{k_1} x) + e(-\l^{(1)}_{k_1} x) \right)^{2l+1}
$.
    Every term contributing to the coefficient of $e(\l^{(1)}_1 x)$ arises in the following way.
    First of all choose $e(\l^{(1)}_1 x)$ from some bracket.
    It can be done in $(2l+1)$ ways.
    Secondly we choose $e(\l^{(1)}_u x)$
    $u$ from some other bracket.
    This  can be done in $2k_1$ ways.
    Clearly, it must be balanced by choosing $e(-\l^{(1)}_u x)$
    from some other bracket.
    There are at most $(2l-1)$ ways of doing this.
    And so on.
    Thus the coefficient of $e(\l^{(1)}_1 x)$ from $j=l$ does not exceed
$$
    \frac{1}{2\sqrt{k_1}} \cdot \frac{1}{(64 k_1)^l l!} \m (2l+1) \m 2k_1 \m (2l-1) \m 2k_1 \m (2l-3) \m \dots \m 2k_1 \m 1
        \le \frac{1}{2\sqrt{k_1}} \frac{2l+1}{2^{4l}} \,.
$$
    Hence
\begin{equation}\label{tmp:(4)}
    \frac{1}{4\sqrt{k_1}} \le \frac{1}{2\sqrt{k_1}} \left( 1 - \sum_{j=1}^{\infty} \frac{2j+1}{2^{4j}} \right)
    \le Q(1,0,\dots,0) \le \frac{1}{2\sqrt{k_1}} \left( 1 + \sum_{j=1}^{\infty} \frac{2j+1}{2^{4j}} \right)
    \le \frac{1}{\sqrt{k_1}}
\end{equation}
 and inequality (\ref{tmp:(2)}) is proved.
    Rather more accurate calculation shows that
    \begin{equation}\label{f:delicate}
        Q(1,0,\dots,0) \le \frac{1}{2\sqrt{k_1}} \,.
    \end{equation}
    Indeed the term  from $j=1$ in (\ref{tmp:(3)}) is negative,
    and its absolute value is equal to
    $$
        \frac{1}{2\sqrt{k_1}} \cdot \frac{1}{64 k_1} (3 (2k_1 - 2) + 3) \ge \frac{1}{2\sqrt{k_1}} \cdot \frac{1}{16} \,.
    $$
    Whence
    $$
        Q(1,0,\dots,0) \le \frac{1}{2\sqrt{k_1}} \left( 1 - \frac{1}{16} + \sum_{j=2}^{\infty} \frac{2j+1}{2^{4j}} \right)
            \le \frac{1}{2\sqrt{k_1}}
    $$
 and inequality (\ref{f:delicate}) is proved.

    It easy to see that there exists $\gamma \in [1/2, 1]$ such that
    $\sum_x f(x) = \d N$, where $f(x) = \gamma g(x)$.
    By assumption $|\L| \le 2^{-12} (\d/\a)^2 \log (1/\d)$.
    Since $f = \gamma g$, it follows that $i\in [p]$, and for any $j\in [k_i]$, we have
    \begin{equation}\label{tmP1}
        |\F{f} (\l^{(i)}_j) | = |\F{f} (-\l^{(i)}_j) | \ge 2\a N
    \end{equation}
    (here we have made use (\ref{tmp:five}) and (\ref{tmp:(4)})).
    Now take $A$ as in Lemma \ref{l:char_f}.
    Then $|A| = [\sum_{x\in \Z_N} f(x)] = [\d N]$, and for all $r\in \Z_N \setminus \{ 0 \}$, we have
    \begin{equation}\label{f:approx_char}
        | \F{A} (r) - \F{f} (r) | \le 20 \sqrt{N} \,.
    \end{equation}
    Using (\ref{tmP1}) and $\a > 40 N^{-1/2}$, we get
    $\{ 0 \} \bigsqcup \L \bigsqcup -\L \subseteq \r_\a (A)$.

    Let $r \notin \{ 0 \} \bigsqcup \L \bigsqcup -\L$.
    Prove that $r\notin \r_\a (A)$.
    As was noted above it suffices to consider residuals $r$
    such that
$
 r = s^{(1)}_1 \l^{(1)}_1 + \dots + s^{(1)}_{k_1} \l^{(1)}_{k_1} +
                    \dots + s^{(p)}_1 \l^{(p)}_1 + \dots + s^{(p)}_{k_p} \l^{(p)}_{k_p}
$, where for all $i\in [p]$ the following inequalities hold
$\sum_{l=1}^{k_i} |s^{(i)}_l| \le 2k_i +1$. Since $r \notin \{ 0
\} \bigsqcup \L \bigsqcup -\L$, it follows that $\sum_{i=1}^p
\sum_{l=1}^{k_i} |s^{(i)}_l| \ge 2$. Let $\sigma_i =
\sum_{l=1}^{k_i} |s^{(i)}_l|$. Then $\sum_{i=1}^p \sigma_i \ge 2$
and either there exists $i\in [p]$ such that  $\sigma_i \ge 2$ or
there exist $i,j \in [p]$, $i\neq j$ such that $\sigma_i, \sigma_j
\ge 1$.

    Let us prove that for all $i\in [p]$,  we have
\begin{displaymath}
        |Q_i (s^{(i)}_1,\dots, s^{(i)}_{k_i})| \le \left\{ \begin{array}{ll}
                        2                         \,,   & \mbox{ if } \sigma_i = 0, \\
                        \frac{1}{2\sqrt{k_i}}      \,,   & \mbox{ if } \sigma_i = 1, \\
                        \frac{2}{k_i \sqrt{k_i}}  \,,   & \mbox{ otherwise. }
                \end{array} \right.
\end{displaymath}
Clearly it suffices to deal with the case $i=1$. Let $\sigma =
\sigma_1$. If  $\sigma = 0$ then we have the upper bound for
$Q_1$. If  $\sigma = 1$ then (\ref{f:delicate}) implies that $Q_1
\le \frac{1}{2\sqrt{k_i}}$. Let $\sigma \ge 2$. Suppose that $j_0$
is the minimal positive integer $j_0 \in [k_1]$ such that $2 j_0 +
1 \ge \sigma$ (if there is not such $j_0$ then $Q_1 = 0$ and we
are done). It is easy to see that it is unnecessary to deal with
all terms in (\ref{tmp:(3)}) such that  $j<j_0$. Suppose that
$\sigma$ is an odd number. Then $\sigma = 2r+1$, $r\ge 1$. The
absolute value of the coefficient from $j=l \ge j_0$ in
(\ref{tmp:(3)}) does not exceed
$$
    \frac{1}{2\sqrt{k_1}} \cdot \frac{1}{(64 k_1)^l l!} \m \frac{(2l+1)!}{|s^{(1)}_1|! \dots |s^{(1)}_{k_1}|! \cdot (2l+1 - \sigma)!}
                \m
$$
\begin{equation}\label{tmp:alef}
    \m
    2k_1 \m (2l+1-\sigma-1) \m 2k_1 \m (2l+1-\sigma-3) \m \dots \m 2k_1 \m 1 := \rho \,.
\end{equation}
Indeed firstly if $s^{(1)}_1 \ge 0$ then choose $|s^{(1)}_1|$
elements $e(\l^{(1)}_1 x)$ from (\ref{tmp:(3)}), if $s^{(1)}_1 <
0$ then choose $|s^{(1)}_1|$ elements $e(-\l^{(1)}_1 x)$. Further
if $s^{(1)}_2 \ge 0$ then choose $|s^{(1)}_2|$ elements
$e(\l^{(1)}_2 x)$ from (\ref{tmp:(3)}) and if $s^{(1)}_2 < 0$ then
choose $|s^{(1)}_2|$ elements $e(-\l^{(1)}_2 x)$. And so on. This
can be done in $(2l+1)! / (|s^{(1)}_1|! \dots |s^{(1)}_{k_1}|!
\cdot (2l+1 - \sigma)!)$ ways. Secondly we choose $e(\l^{(1)}_u
x)$ from some other bracket. This  can be done in $2k_1$ ways.
Clearly, it must be balanced by choosing $e(-\l^{(1)}_u x)$ from
some other bracket. There are at most $(2l+1 - \sigma - 1)$  ways
of doing this. And so one. Finally we have the inequality
(\ref{tmp:alef}). Note that if $\sigma$ is an even number, $\sigma
\ge 2$ then the number $Q_1 (s^{(1)}_1,\dots, s^{(1)}_{k_1})$
equals zero. Using $\sigma \le 2l+1$, we get
$$
    \rho \le \frac{1}{2\sqrt{k_1}} \cdot \frac{k_1^{l-r}}{k_1^l 2^{5l}}
    \cdot
    \frac{(2l+1)(2l)(2l-1) \dots (2l+1 - \sigma +1)(2l+1 - \sigma - 1) (2l+1 - \sigma -3) \dots 1}{l!}
$$
$$
    \le
    \frac{1}{2\sqrt{k_1}} k_1^{-r} (2l)^{r} \frac{2l+1}{2^{4l}}
    \le
    \frac{1}{2k_1\sqrt{k_1}} \frac{(2l+1)2l}{2^{4l}} \,.
$$
Hence
$$
    |Q_1 (s^{(1)}_1,\dots, s^{(1)}_{k_1})| \le \frac{1}{k_1\sqrt{k_1}} \sum_{l=1}^\infty \frac{(2l+1)2l}{2^{4l}}
        \le \frac{2}{k_1\sqrt{k_1}}
$$
and we obtain the required upper bound for $Q_1$.

If there exists $i\in [p]$ such that  $\sigma_i \ge 2$ then
we get
$|\F{g}(r)| \le 2^{-p} N / (k')^{3/2}$, where $k' = \min \{ k_i
~:~ i\in [p] \}$. Combining (\ref{tmp:five}), (\ref{tmp:low}) and
(\ref{tmp:(4)}), we obtain $
    |\F{g} (r)| \le \a N/4
$. If there exist three $\sigma_i=1$, then using (\ref{tmp:five}),
(\ref{tmp:low}) and (\ref{tmp:(4)}) again, we get $
    |\F{g} (r)| \le \a N/4
$. Finally, let there are exactly two $\sigma_i =1$. Combining
(\ref{tmp:five}), (\ref{tmp:low}) and (\ref{f:delicate}), we have
$$
    |\F{g}(r)| \le \frac{2^{-p} N}{16 k'} \le \frac{\d pN}{4|\L|}  < \frac{3\a N}{4} \,.
$$
Anyway for all $r \notin \{ 0 \} \bigsqcup \L \bigsqcup -\L$, we
get $
    |\F{f} (r)| \le |\F{g} (r)| < 3\a N/4
$. Using (\ref{f:approx_char}), we obtain $
    |\F{A} (r)| < \a N
$ for any $r \notin \{ 0 \} \bigsqcup \L \bigsqcup -\L$. Hence
$\r_\a (A) = \{ 0 \} \bigsqcup \L \bigsqcup -\L$. This completes
the proof.

The following result shows that our Theorem \ref{t:main} is sharp.

\Th \label{t:T2_low_Z_N} {\it
    Suppose that $\d,\a \in (0,1]$ are real numbers , $N$ is a prime number,
    $k$ is a positive integer, $2\le k \le 2^{-1} \log (1/\d)$,
    $32 \d^2 \le \a \le \d/4$
     and
    \begin{equation}\label{cond:N,d,a}
        2k \max \left\{ \log \left( \frac{2^6\d k}{\a^2} \right),\, \log \left( \frac{2^6\d^2}{\a^3} \right) \right\} \le \log N \,.
    \end{equation}
    Then there exists a set $A\subseteq \Z_N$ such that  $\d N \le |A| \le 3\d N$,
    $|\r_\a (A)| \ge \frac{\d}{64\a^2}$ and
    for all $k$, satisfying (\ref{cond:N,d,a}), we have
    $T_k (\r_\a (A)) \le \frac{2^{14 k} \d}{\a^{2k}}$.
}

\Note We prove in Theorem \ref{t:T2_low_Z_N}
that $|\r_\a (A)| \ge \frac{\d}{64\a^2}$. Certainly this lower
bound for $|\r_\a (A)|$ is absolutely indispensable because
Theorem \ref{t:T2_low_Z_N} is trivial otherwise.

To prove Theorem \ref{t:T2_low_Z_N} we need in the following
definition and lemma.

\Def \label{d:tL_k,s}
    Let $k,s$ be positive integers.
   Consider the family
   $\t{\L}(k,s)$ of subsets of $\Z_N$.
   A set $\t{\L} = \{ \t{\l}_1, \dots, \t{\l}_{|\t{\L}|} \}$ belongs to the family
   $\t{\L} (k,s)$ if the equality
   \begin{equation}\label{f:d_tL_k,s}
        \sum_{i=1}^{|\t{\L}|} \t{\l}_i s_i = 0 \pmod N \,, \quad \t{\l}_i \in \L \,, \quad s_i \in \Z \,, \quad
        |s_i| \le s\,, \mbox{ the number of } s_i \neq 0 \mbox{ at most } k \,,
   \end{equation}
   implies that all $s_i$ are equal to zero.

   Clearly $\L (ks,s) \subseteq \t{\L} (k,s) \subseteq \L (k,s)$.

\Lemma \label{l:diss_set_exists} {\it
    Let $N,t,k,s$ be positive integers, $k \le t$ and $N > \binom{t}{k} (2s+1)^k$.
    Then $\t{\L} (k,s)$ contains a set of the cardinality $t$.
}
\\
\Proof
    Consider all tuples $(a_1,\dots, a_t)$, where $a_i \in \Z_N$.
    Obviously that there exist $N^t$ of this tuples.
    Further there are at most $\binom{t}{k} (2s+1)^k$ of equations (\ref{f:d_tL_k,s})
    with coefficients $s_1,\dots,s_t$.
    The number of solutions of any non--trivial equation (\ref{f:d_tL_k,s})
    does not exceed $N^{k-1} N^{t-k} = N^{t-1}$.
    Besides
    $$
        N^{t-1} \binom{t}{k} (2s+1)^k < N^t \,.
    $$
    Whence there exists a tuple $(a_1,\dots,a_t)$,
    satisfying the trivial equation only.
    It is easy to see that all residuals in $(a_1,\dots,a_t)$ are different.
    Hence the set $\t{\L} = \{ a_1,\dots,a_t \}$
    belongs to the family $\t{\L} (k,s)$.
    This completes the proof of the lemma.

{\bf Proof of Theorem \ref{t:T2_low_Z_N}.}
    Let $k_1 = 2k$, $t= [ \d/\a ]$, $\eps = \d /t$, $m= \max \{ t,k_1 \}$, $s = \lceil 8m /\eps \rceil$.
    By assumption $2k \log (\frac{2^6\d^2}{\a^3}) \le \log N$ and $2k \log(\frac{2^6\d k}{\a^2}) \le \log N$.
    Hence $N> \binom{t}{k_1} (2s+1)^{k_1}$.
    Using Lemma \ref{l:diss_set_exists}, we find a set $\L = \{ \l_1, \dots, \l_t \}$
    such that $\L$ belongs to the family $\t{\L} (k_1,s)$.

    For any $\l \in \Z_N$ consider one--dimensional Bohr set
    \begin{equation}\label{f:def_Bohr}
        B_\l = B_\l (\eps) = \{ x\in \Z_N ~:~ \left\| \frac{x\l}{N} \right\| \le \eps \}\,,
    \end{equation}
    (see \cite{Bu} for example).
    Clearly,
    $B_\l (\eps) = \{ 0, \pm \l^{-1}, \dots, \pm [\eps N] \l^{-1} \}$.
    Hence  $|B_\l (\eps)| = 2 [\eps N] + 1$.
    By $B_\l^{s} = B_\l^{s} (\eps)$ denote the set $B_\l^{s} = B_\l + s$.
    We shall construct a family of sets $B^{s_1}_{\l_1}, \dots, B^{s_t}_{\l_t}$,
    where $\l_i \in \L$, $s_i \in \Z_N$.
    Let $s_1 =0$ and we obtain the set $B^{s_1}_{\l_1}$.
    Suppose that we have the sets $B^{s_1}_{\l_1}, \dots, B^{s_d}_{\l_d}$.
    Let us construct a residual $s_{d+1}$ and a set $B^{s_{d+1}}_{\l_{d+1}}$.
    Let $C_d = \bigcup_{i=1}^d B^{s_i}_{\l_i}$.
    Clearly, $|C_d| \le d (2[\eps N] + 1) \le t (2 [\eps N] + 1) \le 3\d N$.
    Let $s_{d+1}$ be a residual such that
    \begin{equation}\label{tmp:23:09_16Aug'}
        |B^{s_{d+1}}_{\l_{d+1}} \bigcap C_{d}| \le (2\eps N + 1)^2 t \le 8 \eps \d N \,.
    \end{equation}
    Since
    $$
        \sum_{s\in \Z_N}  |C_d \cap B^{s}_{\l_{d+1}} | = |C_d| |B_s| \,,
    $$
    it follows that such a residual $s_{d+1}$ exists.
    So we have the sets $B^{s_1}_{\l_1}, \dots, B^{s_t}_{\l_t}$.
    Let $A = C_t = \bigcup_{i=1}^t B^{s_i}_{\l_i}$.
    Clearly, $|A| \le 3\d N$.
    Prove that $|A| \ge \d N$.
    Using (\ref{tmp:23:09_16Aug'}), we get
    \begin{equation}\label{f:expl_1'}
        |A| = |C_t| = |C_{t-1}| + |B^{s_t}_{\l_t}| - |C_{t-1} \cap B^{s_t}_{\l_t}| \ge |C_{t-1}| + (2[\eps N] + 1) - 8 \eps \d N
        \ge
    \end{equation}
    \begin{equation}\label{f:expl_2'}
        \ge
            |C_{t-2}| + 2 (2[\eps N] + 1) - 2 \cdot 8 \eps \d N \ge \dots \ge t (2[\eps N] + 1) - t 8 \eps \d N
            \ge  t \eps N  = \d N\,.
    \end{equation}

    Let us prove that $|\r_\a (A)| \ge \frac{\d}{64\a^2}$.
    Let $a \in \Z_N$.
    Assuming that $a$ belongs to reduced residue system, we denote by $|a|$
    the absolute value of $a$.
    We have $|a| \le N/2$ for all $a \in \Z_N$.
    Let $r \in \Z_N$, $r\neq 0$.
    Using the inequality $|1-e^{i \theta}| \ge 2|\theta| / \pi$, $\theta \in [-\pi, \pi]$, we get
    \begin{equation}\label{f:odin}
        | \F{B}_\l (r) | = \left| \sum_{l=-[\eps N]}^{[\eps N]} e(\l^{-1} lr) \right|
            = \left| 2 \frac{e(([\eps N] + 1) \l^{-1}r) - 1}{e(\l^{-1}r) - 1}- 1 \right|
                \le
                    \frac{4}{|e(\l^{-1}r) - 1|} \le \frac{N}{|\l^{-1}r|} \,.
    \end{equation}
    Let us obtain a lower bound for $\F{B}_\l (r)$.
    Let $\l$ belongs to $\L$ and let
    \begin{equation}\label{}
        M_\l = \{\, x\in \Z_N ~:~ x = \l p\,, \quad |p| \le \frac{1}{16\eps} \,\} \,.
    \end{equation}
    Observe that $|M_\l| = 2 [1/(16\eps)] + 1$.
    For all $r\in M_\l$, we have
    \begin{equation}\label{f:dva}
        \F{B}_\l (r) = 2 \sum_{l=0}^{[\eps N]} \cos (2\pi \l^{-1} r l /N) - 1 \ge 2 ([\eps N] + 1) - 1 - ([\eps N] + 1)/4
        \ge \frac{3}{2} \eps N \,.
    \end{equation}
    Formulas (\ref{f:odin}) and (\ref{f:dva}) can be used to
    calculate Fourier coefficients of the sets
    $B_\l = B_\l (\eps)$.
    Note that for any $s,r \in Z_N$, we have $|\F{B}_\l^s (r)| = |\F{B}_\l (r)|$.
    Thus (\ref{f:odin}), (\ref{f:dva}) can be used to find
    the absolute values of Fourier coefficients of the sets $B_\l^s$ too.

    It is easy to see that for all $i,j \in [t]$, $i\neq j$, we have
    $M_{\l_i} \cap M_{\l_j} = \{ 0 \}$.
    Indeed, by assumption the set $\L$ belongs to the family $\t{\L} (k_1,s)$
    and $s\ge 1/ (16\eps)$.
    Hence the only solution of the equation $\l_i p_i = \l_j p_j$, $i\neq j$, $|p_i|, |p_j| \le 1/(16\eps)$
    is $p_i = p_j = 0$.
    Prove that $\bigcup_{i=1}^t M_{\l_i} \subseteq \r_\a (A)$.
    Clearly, $0 \in \r_\a (A)$.
    Let  $i\in [t]$ be an arbitrary number, and $r$ be a non--zero residual
    such that $r$ belongs to some $M_{\l_i}$.
    We have
    $$
        \F{A} (r) = \F{B}_{\l_t}^{s_t} (r) + \F{C}_{t-1} (r) + 8 \theta \eps \d N \,,
    $$
    where $|\theta| \le 1$.
    By the same arguments as in (\ref{f:expl_1'}) --- (\ref{f:expl_2'}), we get
    \begin{equation}\label{f:calculation_A_hat}
        \F{A} (r) = \sum_{l=1}^t \F{B}_{\l_l}^{s_l} (r) + 8 \t{\theta} \eps \d t N \,,
    \end{equation}
    where $|\t{\theta}| \le 1$.
    We have $r\in M_{\l_i}$.
    Using (\ref{f:dva}), we obtain
    $|\F{B}_{\l_i}^{s_i} (r)| \ge 3 \eps N/2$.
    Let $r = \l_i p_i$, $|p_i| \le 1/(16\eps)$, and
    $j\in [t]$ be a number, $j\neq i$.
    Let $p := \l_j^{-1} r = \l_j^{-1} \l_i p_i$.
    Then $\l_j p = \l_i p_i$.
    Since the set $\L$ belongs to the family $\t{\L} (k_1,s)$, it follows that $|p| > s$.
    Using the last inequality and (\ref{f:odin}), we obtain
    \begin{equation}\label{tmp:awful_tmp}
        |\F{B}_{\l_j}^{s_j} (r)| \le \frac{N}{s}\,, \quad r\in M_{\l_i}, \quad r\neq 0 \,.
    \end{equation}
    Combining (\ref{tmp:awful_tmp}) and (\ref{f:calculation_A_hat}), we get
    $$
        |\F{A} (r)| \ge \frac{3}{2} \eps N - \sum_{j\neq i} |\F{B}_{\l_j}^{s_j} (r)| - 8 \eps \d t N
        \ge \frac{3}{2} \eps N - \frac{tN}{s} - 8 \eps \d t N \ge \a N \,.
    $$
    Hence $\bigcup_{i=1}^t M_{\l_i} \subseteq \r_\a (A)$ and
    $|\r_\a (A)| \ge \sum_{i=1}^t |M_{\l_i}| - t = 2t [1/(16\eps)] \ge \frac{\d}{64\a^2}$.

    Finally, we shall show that for all $k$, $2\le k \le 2^{-1} \log (1/\d)$,
    we have
    $T_k (\r_\a (A)) \le \frac{2^{14 k} \d}{\a^{2k}}$.
    Let $g$ be a real number,
    $\l$ belongs to the set $\L$ and let
    \begin{equation}\label{}
        L_\l (g)= \{\, x\in \Z_N ~:~ x = \l p\,, \quad |p| \le g \,\}\,.
    \end{equation}
    Let also $M'_\l = L_\l (8/\eps)$.
    Then $|M'_\l| \le 32/\eps$.
    Let us prove that $\r_\a (A) \subseteq \bigcup_{\l \in \L} M'_\l$.
    Assume the converse.
    Let $r\in \r_\a (A) \setminus \{ 0 \}$ and $r\notin \bigcup_{\l \in \L} M'_\l$.
    If $r\notin \bigcup_{\l \in \L} L_\l (s)$ then  using (\ref{f:calculation_A_hat}), we get
    $$
        |\F{A} (r)| \le \frac{tN}{s} + 8\eps \d t N \le \frac{\eps}{2} N < \a N
    $$
    and $r\notin \r_\a (A)$.
    Let now $r \in \bigcup_{\l \in \L} L_\l (s)$.
    We have $\L \in \t{\L} (k_1,s)$.
    Using this fact, we obtain that for all $\l_1, \l_2 \in \L$, $\l_1 \neq \l_2$
    the following holds
    $L_{\l_1} (s) \cap L_{\l_2} (s) = \{ 0 \}$.
    Let $r\neq 0$.
    It is easy to see that there exists the only $i\in [t]$ such that $r\in L_{\l_i} (s)$.
    Using (\ref{f:calculation_A_hat}), we obtain
    \begin{equation}\label{f:tri}
        |\F{A} (r)| \le |\F{B}_{\l_i}^{s_i} (r)| + \frac{tN}{s} + 8\eps \d t N \,.
    \end{equation}
    We have $r\notin \bigcup_{\l \in \L} M'_\l$.
    Using (\ref{f:odin}), we get
    $|\F{B}_{\l_i}^{s_i} (r)| \le N/g \le \eps N /8$.
    Substituting the last inequality in (\ref{f:tri}), we obtain
    $|\F{A} (r)| \le \eps N/2 < \a N$ and $r\notin \r_\a (A)$.
    Whence $\r_\a (A) \subseteq \bigcup_{\l \in \L} M'_\l$.

    Consider the equation
    \begin{equation}\label{f:TeN''}
        r_1 + \dots + r_k = r'_1 + \dots + r'_k \,,
    \end{equation}
    where all $r_j$, $r'_j$ belong to $\r_\a (A)$.
    We have $\r_\a (A) \subseteq \bigcup_{i=1}^t M'_{\l_i}$.
    Hence any residual from (\ref{f:TeN''}) belongs to a set $M'_{\l_i}$.

    Let $z$ be a non--negative integer,
    and
    $s_1, \dots, s_l$ be positive integers such that
    $s_1 + \dots + s_l + z = 2k$.
    Recall that for all $\l_1, \l_2 \in \L$, $\l_1 \neq \l_2$, we have
    $L_{\l_1} (s) \cap L_{\l_2} (s) = \{ 0 \}$.
    Hence for any $i,j \in [t]$, $i \neq j$, we get
    $M'_{\l_i} \cap M'_{\l_j} = \{ 0 \}$.
    Let $M'_i = M'_{\l_i}$, $i\in [t]$, and $w=2^5 /\eps$.
    Then for all $i\in [t]$, we have $|M'_i| \le w$.
    By $E(s_1,\dots,s_l,z)$ denote the set of all solutions
    $r_1,\dots, r_k$, $r'_1,\dots, r'_k$ of (\ref{f:TeN''})
    such that  among $r_j$, $r'_j$
    there exist exactly  $z$ of zeroes,
    there exist exactly $s_1$ non--zero residuals belong to a set $M'_{j_1}$,
    there exist exactly $s_2$ non--zero residuals belong to a set $M'_{j_2}$,
    $\dots$,
    there exist exactly  $s_l$ non--zero residuals belong to a set $M'_{j_l}$
    and at the same time all sets
    $M'_{j_1}, M'_{j_2}, \dots, M'_{j_l}$ are different.
    Using $\L \in \t{\L} (k_1,s)$, we obtain
    $$
        T_k (\r_\a (A)) = \sum_{l=1}^{2k} \sum_{z=0}^{2k}~ \sum_{s_1,\dots,s_l,\,\, s_1+\dots+s_l+z = 2k}
                                |E(s_1,\dots,s_l,z)|
                        \le
    $$
    \begin{equation}\label{f:-10''}
                        \le
                                tw^{2k-1}
                                    +
                                        \sum_{l=2}^{2k} \sum_{z=0}^{2k}~ \sum_{s_1,\dots,s_l,\,\, s_1+\dots+s_l+z = 2k}
                                            |E(s_1,\dots,s_l,z)| \,.
    \end{equation}
    Let us fixed $s_1, \dots, s_l,z$ and consider the solutions of (\ref{f:TeN''})
    belong to fixed subsets $M'_{j_1}, M'_{j_2}, \dots, M'_{j_l}$.
    Denote by $E(s_1,\dots,s_l,z) (M'_{j_1}, M'_{j_2}, \dots, M'_{j_l})$
    the set of
    all these solutions.
    Rewrite the equation (\ref{f:TeN''}) as
    \begin{equation}\label{}
        u_1 + \dots + u_l = 0 \,,
    \end{equation}
    where $u_i \in M'_{j_i}$, $i\in [l]$.
    Since $\L$ belongs to $\t{\L} (k_1,s)$,
    it follows that all residuals $u_i$ equal zero.
    Hence, we have
    $$
        | E(s_1,\dots,s_l,z) (M'_{j_1}, M'_{j_2}, \dots, M'_{j_l}) |
            \le
                \frac{(2k)!}{s_1! \dots s_l! z!} w^{s_1-1} \m \dots \m w^{s_l-1}
                    \le
                         \frac{(2k)!}{s_1! \dots s_l! z!} w^{2k - l} \,.
    $$
    Whence
    \begin{equation}\label{f:-11''}
        |E(s_1,\dots,s_l,z)|
            \le
                \binom{t}{l} \frac{(2k)!}{s_1! \dots s_l! z!} w^{2k - l}
                    \le
                        \frac{t^l}{l!} \cdot \frac{(2k)!}{s_1! \dots s_l! z!} w^{2k - l} \,.
    \end{equation}
    Combining (\ref{f:-11''}) and (\ref{f:-10''}), we get
    $$
        T_k (\r_\a (A))
            \le
                tw^{2k-1}
                    +
                        \sum_{l=2}^{2k} \frac{t^l}{l!} w^{2k - l}
                            \sum_{z=0}^{2k}~ \sum_{s_1,\dots,s_l,\,\, s_1+\dots+s_l+z = 2k} \frac{(2k)!}{s_1! \dots s_l! z!}
            \le
    $$
    \begin{equation}\label{}
            \le
                tw^{2k-1}
                    +
                        \sum_{l=2}^{2k} \frac{t^l}{l!} w^{2k - l} (l+1)^{2k}
            =
                tw^{2k-1}
                    +
                        w^{2k} \sum_{l=2}^{2k} \left( \frac{t}{w} \right)^l \cdot (l+1)^{2k} \cdot \frac{1}{l!} \,.
    \end{equation}
    Consider the function $f(l) = (t/w)^l (l+1)^{2k}$.
    It is easy to see that $f(l)$ has maximum at $l_0 = 2k / \ln (w / t) - 1$
    and for all $l\ge l_0$ the function $f(l)$ is monotonically decreasing.
    By assumption $k \le 2^{-1} \log (1/\d)$.
    Hence $l_0 \le 1$.
    It follows that
    $$
        T_k (\r_\a (A))
            \le
                t w^{2k-1} + 2^{2k} t w^{2k-1}
                    \le
                        2^{2k+1} t w^{2k-1}
                            \le
                               \frac{2^{14 k} \d}{\a^{2k}} \,.
    $$
    This completes the proof of the Theorem.

\refstepcounter{section} \label{positive}

{\bf \arabic{section}. Proof of Theorem \ref{t:Chang_log_m}.}

  We assume in the section that $N$ is a prime number.

\Def \label{def:L_d} {
    Let $k,s,d$ be positive integers.
    Let also $\L = \{ \l_1, \dots, \l_{|\L|} \} \subseteq \Z_N$ be a set such that
    $\L \cap -\L = \emptyset$.
    Let
    $\v{v}_1 = ( v_1^{(1)}, \dots, v_1^{(d)}), \dots, \v{v}_{|\L|} = ( v_{|\L|}^{(1)}, \dots, v_{|\L|}^{(d)})$
    be vectors from $\Z^d$ such that for all $i\in [d]$, $j\in [|\L|]$, we have $|\v{v}_{j}^{(i)}| \le s$.
    Consider the equation
    \begin{equation}\label{f:def1}
        \l_1 \v{v}_1 + \dots +\l_{|\L|} \v{v}_{|\L|} = 0 \pmod N \,,
    \end{equation}
    where $\l_i \in \L$
    and
    for any $i\in [d]$, we have $\sum_{j=1}^{|\L|} |v_j^{(i)} | \le k$.
    Consider the family $\L_d (k,s)$ of subsets of $\Z_N$.
    Our set $\L$ belongs to the family $\L_d (k,s)$, if any equation (\ref{f:def1})
    imply that the matrix
    \begin{displaymath}
    \left( \begin{array}{ccc}
        v_1^{(1)} & \cdots & v_{|\L|}^{(1)} \\
        \cdots    & \cdots & \cdots \\
        v_1^{(d)} & \cdots & v_{|\L|}^{(d)}
    \end{array} \right)
    \end{displaymath}
    has the rank at most $d-1$.
}

As was noted above the definition of $\L_1 (k,1)$ can be found in
\cite{Ruzsa_independent}, and the definition of  $\L_1 (k,s)$ can
be found in \cite{Sh_exp1}.

    For an arbitrary $\L \in \L_d (k,s)$, we obtain the following upper bound for $T_k (\L)$.

\St
    \label{st:sol}
    {\it
        Let $N, k, s, d$ be positive integers, $s\ge 3$,
        $N\ge s+1$ be a prime number,
 and
        $\L \subseteq \Z_N$ be a subset of the family
        $\L_d (2k,s)$.
        Then
        \begin{equation}\label{f:T_k_up_estimate}
            T_k (\L) \le 2^{9k} k^k |\L|^k (s+1)^{2d} \cdot
                2^{ \frac{2s k (\log k)^2 }{\log (k^{2s} |\L|^{s-2})}} \,.
        \end{equation}
    }

    \Exm Let $k\ge 2$, $\log |\L| \ge \log^2 k$, and $\L$ be an arbitrary subset of the family $\L_d (k,3)$.
    Using  (\ref{f:T_k_up_estimate}), we get
    $T_k (\L) \le 2^{20k + 4d} k^k |\L|^k$.

    {\bf Proof of Statement \ref{st:sol}. }
    Let $x\in \Z_N$ be a residual.
    By $N_k (x)$ define the number of vectors $(\l_1,\dots, \l_k)$ such that
    all $\l_i$
    belong to $\L$ and
    \begin{equation}\label{f:sum_l=x}
        \l_1 + \dots + \l_k = x \,.
    \end{equation}
    Then $T_k (\L) = \sum_{x\in \Z_N} N^2_k (x)$.
    Let $s_1,\dots, s_l$ be positive integers such that $s_1 + \dots + s_l = k$.
    By $E(s_1,\dots,s_l) (x)$ denote the set of all solutions $(\l_1,\dots,\l_k)$
    of (\ref{f:sum_l=x})
     such that
     among $\l_1,\dots,\l_k$ there exist exactly $s_1$ residuals equal $\t{\l}_1$,
    there exist exactly $s_2$  residuals equal $\t{\l}_2$, $\dots$,
    there exist exactly $s_l$  residuals equal $\t{\l}_l$
    such that
    $s_1 \t{\l}_1 + \dots + s_l \t{\l}_l = x$
    and all $\t{\l}_i$ are different.
    Let us denote the set $E(s_1,\dots,s_l) (x)$ by $E({\v{s}}) (x)$ for simplicity.
    Recall that for $s_1,\dots,s_l$ in the definition of
    $E({\v{s}}) (x) = E(s_1,\dots,s_l) (x)$ the following equality holds : $\sum_{i=1}^l s_i = k$.
    We have
    $$
        N_k (x) = \sum_{\v{s}} |E({\v{s}}) (x)| \,.
    $$
    Whence
    \begin{equation}\label{}
        \sigma = T_k (\L) = \sum_{x\in \Z_N} ( \sum_{\v{s}} |E({\v{s}}) (x)| )^2 \,.
    \end{equation}
    Let $\v{s} = (s_1, \dots, s_l)$
    and
    $G = G(\v{s}) = \{ i ~:~ s_i \le s \}$,
    $B = B(\v{s}) = \{ i ~:~ s_i > s \}$.
    Then $|G(\v{s})| + |B(\v{s})| = l(\v{s}) = l$.
    We have
    \begin{equation}\label{f:zero}
        l \le k - s |B| \,.
    \end{equation}
    Indeed
    \begin{equation}\label{tmp:5:42}
        k = \sum_{i\in G} s_i + \sum_{i\in B} s_i \ge |G| + (s+1) |B| = l + s |B| \,.
    \end{equation}
    Using (\ref{tmp:5:42}), we obtain (\ref{f:zero}).
    Let also
    $$
        l_j = l_j (\v{s}) = |\, \{ i ~:~ s_i = j,~ i\in [l] \} \, | \,, \quad \quad \quad j = 1,2,\dots, r=r(\v{s})\,, \quad r\neq 0  \,.
    $$

    The next lemma was proved in \cite{Sh_exp1}.

    \Lemma
    \label{l:point2}
    {\it
            For all $\v{s}$, $\sum_{i=1}^l s_i = k$ the number of $x\in \Z_N$ such that
            $E({\v{s}}) (x) \neq \emptyset$
            does not exceed $|\L|^l / l_1!$.
    }

    We need in two lemmas.

    \Lemma
    \label{l:basis}
    {\it
        Suppose that $n,t,s$ are positive integers, $t\le n$,
        $\v{u}_1, \dots, \v{u}_t \in \Z_N^n$ are linearly--independent vectors over $\Z_N$, and $N\ge s+1$.
        Let
        $$
            Q(s) := \{ \v{x} = (x_1, \dots, x_n) \in \Z_N^n ~:~
                                x_i \in \{ 0,1,\dots, s\},\,\,\, i=1,\dots,n \}
        $$
        be a $n$--dimensional cube
         and
        $L = \{ \v{x} \in \Z_N^n ~:~ \v{x} = \sum_{i=1}^t m_i \v{u}_i,\,\,\, m_i \in \Z_N \}$.
        Then
            $
                |L \cap Q(s)| \le (s+1)^t
            $.
    }
    \\
    {\bf Proof of Lemma \ref{l:basis}. }
    Let $\v{u}_1 = ( u_1^{(1)}, \dots, u_1^{(n)} ), \dots, \v{u}_t = ( u_{t}^{(1)}, \dots, u_{t}^{(n)} )$.
    Note that the cube $Q(s)$ is invariant under permutations of coordinates.
    So we can assume without loss of generality that
    the vectors $\v{u}_1, \dots, \v{u}_t$ have the form
    $\v{u}_1 = ( 1,\dots, 0,0, u_{t}^{(t+1)}, \dots, u_1^{(n)} ),
     \v{u}_2 = ( 0,1,\dots, 0, u_{t}^{(t+1)}, \dots, u_2^{(n)} ),
     \dots,
     \v{u}_t = ( 0,\dots, 0,1, u_{t}^{(t+1)}, \dots, u_{t}^{(n)} )$.
    Let $\v{x}$ be an arbitrary vector, $\v{x} \in L \cap Q(s)$.
    Then there exist residuals $m_1, \dots, m_t$ such that
    $\v{x} = \sum_{i=1}^t m_i \v{u}_i$.
    Clearly, $m_1 = x_1, \dots, m_t = x_t$.
    Since $\v{x} \in Q(s)$, it follows that $m_i \in \{0,1,\dots,s\}$, $i=1,\dots,t$.
    Hence $|L \cap Q(s)| \le (s+1)^t$ as required.

    \Lemma
    \label{l:sol}
    {\it
          For any $\v{s}$, $\sum_{i=1}^l s_i = k$ and for any $x\in \Z_N$, we have
            \begin{equation}\label{}
                |E({\v{s}}) (x)| \le \frac{k!}{s_1! \dots s_l!} (s+1)^d |\L|^{|B(\v{s})|} \,.
            \end{equation}
    }
    \\
    {\bf Proof of Lemma \ref{l:sol}. }
    Let $(\l^{(1)},\dots, \l^{(k)}) \in E({\v{s}}) (x)$.
    Then $\sum_{i=1}^k \l^{(i)} = \sum_{i=1}^l s_i \t{\l}^{(i)} = x$
    and all residuals $\t{\l}^{(i)} \in \{ \l^{(1)},\dots, \l^{(k)} \}$ are different.
    Thus for any $(\l^{(1)},\dots, \l^{(k)}) \in E({\v{s}}) (x)$
    there exists the only vector
    $(\t{\l}^{(1)}, \dots, \t{\l}^{(l)})$ such that all $\t{\l}^{(i)}$ are different.
    Let us fixed $\t{\l}^{(i)}$, $i\in B(\v{s})$.
    We obtain the set $K = K (B(\v{s}))$ of $(\l^{(1)},\dots, \l^{(k)}) \in E({\v{s}}) (x)$.
    Let us prove that the cardinality of  $K$ does not exceed $(s+1)^d \frac{k!}{s_1! \dots s_l!}$.

    Let $(\l^{(1)},\dots, \l^{(k)})$ be a vector from $K$.
    We have
    \begin{equation}\label{f:pos_1}
        \sum_{i\in G(\v{s})} s_i \t{\l}^{(i)} = x - \sum_{i\in B(\v{s})} s_i \t{\l}^{(i)} = x' \,.
    \end{equation}
    Since the elements $\t{\l}^{(i)}$, $i \in B(\v{s})$ are fixed, it follows that the residual $x'$
    is the same for all $(\l^{(1)},\dots, \l^{(k)}) \in K$.

    Let the set $\L = \{ \l_1, \dots, \l_{|\L|} \}$
    be ordered in an arbitrary way.
    To each $(\l^{(1)},\dots, \l^{(k)}) \in K$ assign the vector
    $\v{u} = (u_1, \dots, u_{|\L|})$, where
    \begin{displaymath}
        u_j = \left\{ \begin{array}{ll}
                s_i \,, & \mbox{ if for some } i\in G(\v{s}) \mbox{ we have } \l_j = \t{\l}^{(i)}, \\
                0 \,,   & \mbox{ otherwise. }
                \end{array} \right.
    \end{displaymath}
    Let $\v{w} = (w_1, \dots, w_m)$, $\v{w}' = (w_1',\dots, w_m')$ be two vectors from $\Z^m$.
    By $(\v{w}, \v{w}')$ denote the inner product of these vectors:
    $(\v{w}, \v{w}') = \sum_{j=1}^m w_j w_j'$.
    Rewrite (\ref{f:pos_1}) as
    \begin{equation}\label{f:second}
        (\v{u}, \v{\l}) = x' \,,
    \end{equation}
    where $\v{\l} = (\l_1, \dots, \l_{|\L|})$.
    Let $\v{u}_1, \dots, \v{u}_t$ be a maximal system of vectors which are linearly--independent over $\Z_N$
    and
    such that any of these vectors satisfies  (\ref{f:second}).

    Suppose that $t\ge d+1$.

    Since all vector $\v{u}_i$, $i=1,\dots, t$ satisfy (\ref{f:second}), it follows that
    \begin{displaymath}
    \left\{ \begin{array}{ll}
        (\v{u}_2 - \v{u}_1, \v{\l}) = 0 \,.\\
        \dots \dots \dots \dots \dots \dots \\
        (\v{u}_{d+1} - \v{u}_1, \v{\l}) = 0 \,.
    \end{array} \right.
    \end{displaymath}
    Rewrite this system as $\l_1 \v{v}_1 + \dots + \l_{|\L|} \v{v}_{|\L|} = 0$,
    where $\v{v}_i$ are vectors from $\Z^d$.
    Let $\v{v}_1 = ( v_1^{(1)}, \dots, v_1^{(d)} ), \dots, \v{v}_{|\L|} = ( v_{|\L|}^{(1)}, \dots, v_{|\L|}^{(d)} )$.
    Since $\L \cap -\L = \emptyset$, it follows that  $|\v{v}^{(j)}_i| \le s$.
    It is easy to see that for any $i\in [d]$, we have $\sum_{j=1}^{|\L|} |v_j^{(i)} | \le 2k$.
    Consider the matrix
    \begin{displaymath}
    M =
    \left( \begin{array}{ccc}
        v_1^{(1)} & \cdots & v_{|\L|}^{(1)} \\
        \cdots    & \cdots & \cdots \\
        v_1^{(d)} & \cdots & v_{|\L|}^{(d)}
    \end{array} \right)
    \end{displaymath}
    Let $\v{p}_j = (v_1^{(j)}, \dots, v_{|\L|}^{(j)})$, $j=1,\dots,d$ be rows of $M$.
    Clearly, $\v{p}_j = \v{u}_{j+1} - \v{u}_1$, $j = 1,\dots, d$.
    Since the vectors $\v{u}_1, \dots, \v{u}_{d+1}$ are linearly--independent, it follows that
    the vectors $\v{p}_1, \dots, \v{p}_d$ are also linearly--independent.
    Hence the rank of $M$ is equal to $d$
    with contradiction to the definition of the family $\L_d (2k,s)$.

    Thus, we have $t \le d$.
    Let $\v{u}$ be an arbitrary vector satisfies (\ref{f:second})
    By maximality of the system $\v{u}_1, \dots, \v{u}_t$, we get
    $\v{u} = \sum_{i=1}^t m_i \v{u}_i$, where $m_i \in \Z_N$.
    All coordinates of $\v{u}$ belong to $\{0,1,\dots, s\}$.
    Using Lemma \ref{l:basis}, we obtain that the number of such $\v{u}$
    does not exceed $(s+1)^d$.
    Clearly, any vector $\v{u}$ corresponds to the tuple $\{ \t{\l}^{(i)} \}_{i\in G(\v{s})}$.
    Besides we fixed residuals $\{ \t{\l}^{(i)} \}_{i\in B(\v{s})}$.
    Using the definition of the set $E(\v{s}) (x)$, we get that the number
    of permutations of the tuple
    $\{ \t{\l}^{(i)} \}_{i\in G(\v{s})} \bigsqcup \{ \t{\l}^{(i)} \}_{i\in B(\v{s})}$
    equals $\frac{k!}{s_1! \dots s_l!}$.
    Hence the cardinality of the set $K$ at most $(s+1)^d \frac{k!}{s_1! \dots s_l!}$
    and the cardinality of $E(\v{s}) (x)$ at most $(s+1)^d \frac{k!}{s_1! \dots s_l!} |\L|^{|B(\v{s})|}$.
    This completes the proof of Lemma \ref{l:sol}.

    Let us return to the proof of Statement \ref{st:sol}.

    Let $t= (k \log k) / \log (k^{2s} |\L|^{s-2})$.
    We can assume without loss of generality that
    \begin{equation}\label{f:T_k_trivial}
        |\L|^k \ge 2^{9k} k^k \,.
    \end{equation}
    Using (\ref{f:T_k_trivial}), we get
$        |\L| \ge 2^9 k \ge k$.
    Let us estimate the sum  $\sigma$.
    \begin{equation}\label{}
        \sigma \le 2 \left( \sum_{x\in \Z_N} \left( \sum_{\v{s} : |B(\v{s})| \le t} |E({\v{s}}) (x)| \right)^2
            +  \sum_{x\in \Z_N} \left( \sum_{\v{s} : |B(\v{s})| > t} |E({\v{s}}) (x)| \right)^2 \right) = 2 \sigma_1 + 2 \sigma_2 \,.
    \end{equation}
    We have
    \begin{equation}\label{}
        \sigma_2 \le \sum_{\v{s}_1, \v{s}_2 : |B(\v{s}_1)| > t, |B(\v{s}_2)| > t} \sum_x |E({\v{s}_1}) (x)| \cdot  |E({\v{s}_2}) (x)|
    \end{equation}
    If $|B(\v{s}_1)| > |B(\v{s}_2)|$ then put $\v{s}^* = \v{s}_1$.
    If $|B(\v{s}_1)| \le |B(\v{s}_2)|$ then set
    $\v{s}^* = \v{s}_2$.
    Let also $P_k (\v{s}) = k! / (s_1! \dots s_l!)$.
    Using Lemma \ref{l:sol}, we obtain $|E({\v{s}_1}) (x)| \le P_k(\v{s}_1) (s+1)^{d} |\L|^{|B(\v{s}^*)|}$
     and $|E({\v{s}_2}) (x)| \le P_k(\v{s}_2) (s+1)^{d} |\L|^{|B(\v{s}^*)|}$.
    Further, using Lemma \ref{l:point2}, we get
    \begin{equation}\label{}
       \sigma_2  \le (s+1)^{2d} \sum_{\v{s}_1, \v{s}_2 : |B(\v{s}_1)| > t, |B(\v{s}_2)| > t} |\L|^{l(\v{s}^*)} |\L|^{2|B(\v{s}^*)|}
                        P_k(\v{s}_1) P_k(\v{s}_2) \,.
    \end{equation}
    Taking into account (\ref{f:zero}), we have
    \begin{equation}\label{}
        \sigma_2 \le (s+1)^{2d} \sum_{\v{s}_1, \v{s}_2 : |B(\v{s}_1)| > t, |B(\v{s}_2)| > t}
                        |\L|^{k-s|B(\v{s}^*)|} |\L|^{2|B(\v{s}^*)|} P_k(\v{s}_1) P_k(\v{s}_2)
                            \le
    \end{equation}
    \begin{equation}\label{}
                            \le (s+1)^{2d} |\L|^k |\L|^{-t(s-2)} \sum_{\v{s}_1, \v{s}_2} P_k(\v{s}_1) P_k(\v{s}_2)
                                \le 2^4 (s+1)^{2d} |\L|^k |\L|^{-t(s-2)} (k^k)^2 \,.
    \end{equation}
    Since
    $t = (k \log k) / \log (k^{2s} |\L|^{s-2})$, it follows that
    \begin{equation}\label{f:k^k_L^-t}
        k^k |\L|^{-t(s-2)} \le 2^{ \frac{2s k (\log k)^2 }{\log (k^{2s} |\L|^{s-2})}} \,.
    \end{equation}
    Hence
    \begin{equation}\label{f:A}
        \sigma_2 \le 2^4 k^k |\L|^k (s+1)^{2d} \cdot 2^{ \frac{2s k (\log k)^2 }{\log (k^{2s} |\L|^{s-2})}} \,.
    \end{equation}
     Let us estimate $\sigma_1$.
    $$
        \sigma_1 \le 2 \left( \sum_{x\in \Z_N} \left( \sum_{\v{s} : |B(\v{s})| \le t, l(\v{s}) \le k - st} |E({\v{s}}) (x)| \right)^2
                                + \sum_{x\in \Z_N} \left( \sum_{\v{s} : |B(\v{s})| \le t, l(\v{s}) > k - st} |E({\v{s}}) (x)| \right)^2 \right)
                                =
    $$
    \begin{equation}\label{}
                                =
                                2\sigma^{'}_1 + 2\sigma^{''}_1 \,.
    \end{equation}
    We have
    \begin{equation}\label{}
        \sigma^{'}_1 \le \sum_{\v{s}_1, \v{s}_2 : |B(\v{s}_1)|, |B(\v{s}_2)| \le t,\,\, l(\v{s}_1), l(\v{s}_2) \le k - st}
                            \sum_x |E({\v{s}_1}) (x)| \cdot  |E({\v{s}_2}) (x)|
    \end{equation}
   Using Lemmas \ref{l:point2}, \ref{l:sol} and (\ref{f:k^k_L^-t}), we obtain
    \begin{equation}\label{}
         \sigma^{'}_1 \le (s+1)^{2d}
                            \sum_{\v{s}_1, \v{s}_2 : |B(\v{s}_1)|, |B(\v{s}_2)| \le t,\,\, l(\v{s}_1), l(\v{s}_2) \le k - st}
                            |\L|^{l(\v{s}^*)} |\L|^{2|B(\v{s}^*)|} P_k(\v{s}_1) P_k(\v{s}_2)
                            \le
    \end{equation}
    \begin{equation}\label{f:B}
                            \le 2^4 (s+1)^{2d} |\L|^k |\L|^{-t(s-2)} (k^k)^2
                                \le 2^4  k^k |\L|^k (s+1)^{2d} \cdot
                                   2^{ \frac{2s k (\log k)^2 }{\log (k^{2s} |\L|^{s-2})}} \,.
    \end{equation}
    We need in an upper bound for $\sigma^{''}_1$.
    For any  $\v{s} = (s_1, \dots, s_l)$,\, $\sum_{i=1}^l s_i = k$, we have
    \begin{equation}\label{tmp:osel_I}
        l_1 + \dots + l_r = l \quad \mbox{ and } \quad l_1 + 2 l_2 + \dots + r l_r = k \,.
    \end{equation}
    Using (\ref{tmp:osel_I}), we get
    $
        l = k  - ( l_2 + 2l_3 + \dots + (r-1) l_r)
    $.
     On the other hand $l\ge k - st$.
    It follows that $l_2 + 2l_3 + \dots + (r-1) l_r \le st$.
    Further,
    $l_2 + \dots + l_r \le l_2 + 2l_3 + \dots + (r-1) l_r \le st$.
    Whence $l_1 = l - (l_2 + \dots + l_r) \ge l - st \ge k - 2st$.
    Taking into account Lemmas \ref{l:point2}, \ref{l:sol} and (\ref{f:zero}), we obtain
    \begin{equation}\label{}
        \sigma^{''}_1 \le \sum_{\v{s}_1, \v{s}_2 : |B(\v{s}_1)|, |B(\v{s}_2)| \le t,\, l(\v{s}_1), l(\v{s}_2) > k - st}
                            \sum_x |E({\v{s}_1}) (x)| \cdot  |E({\v{s}_2}) (x)|
                            \le
    \end{equation}
    \begin{equation}\label{}
                            \le
                                (s+1)^{2d}
                                \sum_{\v{s}_1, \v{s}_2 : |B(\v{s}_1)|, |B(\v{s}_2)| \le t,\,\, l(\v{s}_1), l(\v{s}_2) > k - st}
                                    \frac{|\L|^{l(\v{s}^*)} |\L|^{2|B(\v{s}^*)|}}{l_1 (\v{s}^*) !} P_k (\v{s}_1) P_k (\v{s}_2)
                            \le
    \end{equation}
    \begin{equation}\label{}
                            \le (s+1)^{2d}
                                \sum_{\v{s}_1, \v{s}_2 : |B(\v{s}_1)|, |B(\v{s}_2)| \le t,\,\, l(\v{s}_1), l(\v{s}_2) > k - st}
                                    \frac{|\L|^{k - s|B(\v{s}^*)|} |\L|^{2|B(\v{s}^*)|}}{l_1 (\v{s}^*) !} P_k (\v{s}_1) P_k (\v{s}_2) \,.
    \end{equation}
    Since $l_1 = l_1 (\v{s}^*) \ge k - 2st$, it follows that
    \begin{equation}\label{}
        \sigma^{''}_1 \le (s+1)^{2d} \frac{|\L|^k}{[k-2st]!} \sum_{\v{s}_1, \v{s}_2} P_k (\v{s}_1) P_k (\v{s}_2)
                            \le 2^4 (s+1)^{2d} \frac{|\L|^k}{[k-2st]!} (k^k)^2 \,.
    \end{equation}
    We have  $|\L| \ge k$.
    Whence $t\le k / (3s - 2)$.
    Using the last inequality, we get
    $$
        [k-2st]! \ge [k-2st]^{[k-2st]} / e^{k} \ge k^{[k-2st]} / (8e)^{k} \,.
    $$
    Since
    $t = (k \log k) / \log (k^{2s} |\L|^{s-2} )$, it follows that
    $k^{2st} \le 2^{ \frac{2s k (\log k)^2}{\log (k^{2s} |\L|^{s-2})}}$.
    Further,
    $$
        [k-2st]! \ge k^k / (2^{5k} 2^{ \frac{2s k (\log k)^2}{\log (k^{2s} |\L|^{s-2})}}) \,.
    $$
    Hence
    \begin{equation}\label{f:C}
        \sigma^{''}_1 \le 2^4 2^{5k} k^k |\L|^k (s+1)^{2d}
            \cdot
                2^{ \frac{2s k (\log k)^2}{\log (k^{2s} |\L|^{s-2})}} \,.
    \end{equation}
    Combining (\ref{f:A}), (\ref{f:B}) and (\ref{f:C}), we finally obtain
    \begin{equation}\label{}
        \sigma = T_k (\L) \le 2^{9k} k^k |\L|^k (s+1)^{2d}
                    \cdot
                        2^{ \frac{2s k (\log k)^2}{\log (k^{2s} |\L|^{s-2})}} \,.
    \end{equation}
    This completes the proof of Statement \ref{st:sol}.


   {\bf Proof of Theorem \ref{t:Chang_log_m}}
   Let $k = [\log (1/\d)]$.
   Since $0\in \r_\a$ and $\r_\a = - \r_\a$, it follows that there exists a set $\r^{(1)}_\a$ such that
   $\r_\a = \r^{(1)}_\a \bigsqcup -\r^{(1)}_\a \bigsqcup \{ 0 \}$
   and $\r^{(1)}_\a \cap -\r^{(1)}_\a = \emptyset$.
   Let $s=3$ and $\L = \{ \l_1, \dots, \l_{|\L|} \}$ be a maximal subset of $\r^{(1)}_\a$
   such that $\L$ belongs to $\L_d (2k,s)$.
   Let $\L^* = (\bigcup_{j=1}^3 j^{-1} \L ) \bigcup ( - \bigcup_{j=1}^3 j^{-1} \L )$.
   Then $|\L^*| \le 8 |\L|$.

   Let us prove that for any $r \in \r^{(1)}_\a$ there exists a vector $\v{u} = (u_1, \dots, u_d)$ and
   there exist vectors
   $\v{v}_1 = (v^{(1)}_1, \dots, u^{(d)}_1), \dots, \v{v}_1 = (v^{(1)}_{|\L|}, \dots, u^{(d)}_{|\L|})$
   such that
   $|u_l| \le s$, $l=1,\dots, d$, $|v_j^{(i)}| \le s$, $i=1,\dots,d$, $j=1,\dots, |\L|$ and
   \begin{equation}\label{tmp:4:47}
        r \v{u} = \sum_{i=1}^{|\L|} \l_i \v{v}_i\,,
   \end{equation}
   where for all $i\in [d]$ the following inequality holds $\sum_{j=1}^{|\L|} |v_j^{(i)}| \le k$ and
   the rank of the matrix
    \begin{displaymath}
    M =
    \left( \begin{array}{ccc}
        v_1^{(1)} & \cdots & v_{|\L|}^{(1)} \\
        \cdots    & \cdots & \cdots \\
        v_1^{(d)} & \cdots & v_{|\L|}^{(d)}
    \end{array} \right)
    \end{displaymath}
   equals $d$.

   If such vectors exist then it is easy to see
   that (\ref{f:r=log+}) holds.
   Indeed since the set $\L$ belongs to $\L_d (2k,s)$, it follows that
   the vector $\v{u}$ does not equal zero.
   Hence it has non--zero component.
   Without loss of generality it can be assumed that the first component of $\v{u}$
   does not equal zero.
   For any $i$ take $i$--th equation of system (\ref{tmp:4:47}) such that $\v{u}_i = 0$
   and add the first equation of (\ref{tmp:4:47}) to this equation.
   We obtain a new system
   \begin{equation}\label{tmp:4:47+}
        r \v{u}' = \sum_{i=1}^{|\L|} \l_i \v{v}'_i\,,
   \end{equation}
   where {\it all} components of a vector $\v{u}'$ do not equal zero,
   for any $i\in [d]$ we have $\sum_{j=1}^{|\L|} |(v'_j)^{(i)}| \le 2k \le 4 \log 1/ \d$
   and a matrix $M' = \{ \v{v}'_1, \dots, \v{v}'_{|\L|} \}$
   has the rank $d$.
   Clearly, it can be assumed that all components of $\v{u}'$ belong to $[s]$.
   For any $i\in [|\L|]$ and for any $j\in [s]$, we have $j^{-1} \l_i \in \L^*$.
   Hence system (\ref{tmp:4:47+}) implies (\ref{f:r=log+}) for all $r\in \r^{(1)}_\a$.
   This obviously implies that equation (\ref{f:r=log+}) holds for all $r\in -\r^{(1)}_\a$.

   Thus let $r$ be an arbitrary element of $\r^{(1)}_\a \setminus \L$.
   Let us consider all equations
   \begin{equation}\label{tmp:4:47++}
        \sum_{i=1}^{|\L|} \t{\l}_i \v{v}_i + r \v{u} = \v{0} \,,
   \end{equation}
   such that $|v_j^{(i)}|, |u^{(i)}| \le s$
   and for all $i\in [d]$, we have
   $\sum_{j=1}^{|\L|} |v_j^{(i)}| + |u^{(i)}| \le k$.
   Consider all matrices
   \begin{displaymath}
    M_1 =
    \left( \begin{array}{cccc}
        v_1^{(1)} & \cdots & v_{|\L|}^{(1)} & u^{(1)}\\
        \cdots    & \cdots & \cdots & \cdots \\
        v_1^{(d)} & \cdots & v_{|\L|}^{(d)} & u^{(d)}
    \end{array} \right)
    \end{displaymath}
    If all these matrices have the rank at most $d-1$
    then
    we obtain a contradiction with maximality of $\L$.
    It follows that there exists an equation (\ref{tmp:4:47++}) such that  the rank of $M_1$
    equals $d$.
    Let $M$ be the $(d\m |\L|)$ matrix composed of first $|\L|$ columns of $M_1$.
    Using (\ref{tmp:4:47++}), we get that the rank of $M$ is also equals $d$.
    As was noted above this implies (\ref{f:r=log+}).

     Let us obtain the bound
     $|\L^*| \le \max( \, 2^{30 + 8d (\log (1/\d))^{-1}} \cdot (\d /\a)^2 \log (1/\d), 2^{ ( \log \log (1/\d) )^2 + 3} \, )$.

    If $\log |\L| < (\log k)^2$ then $|\L| \le 2^{ ( \log \log (1/\d) )^2}$
    and $|\L^*| \le 2^{ ( \log \log (1/\d) )^2 + 3}$.
    Suppose that $\log |\L| \ge (\log k)^2$.
    Using Statement \ref{st:sol}, we get
    $T_k (\L) \le 2^{20k + 4d} k^k |\L|^k$.
    On the other hand, using Theorem \ref{t:main}, we obtain
    $T_k (\L) \ge \d \a^{2k} |\L|^{2k} / (2^{4k} \d^{2k})$.
    Hence $|\L| \le 2^{27 + 8d (\log (1/\d))^{-1} } (\d/\a)^2 \log (1/\d)$
    and
    $|\L^*| \le 2^{30 + 8d (\log (1/\d))^{-1}} \cdot (\d/\a)^2 \log (1/\d)$.

    In any case, we have $|\L^*| \le \max( \, 2^{30 + 8d (\log (1/\d))^{-1}} \cdot (\d /\a)^2 \log (1/\d), 2^{ ( \log \log (1/\d) )^2 + 3} \, )$.

    Let us prove that $|\L^*| \le 2^{30} (\d/\a)^2 \log^{\_phi} (1/\d)$.
    If $|\L| < k^\_phi$ then $|\L^*| \le 8 |\L| \le 8 k^\_phi$ and we are done.
    If $|\L| \ge k^\_phi$ then using Statement \ref{st:sol}, we obtain
    \begin{equation}\label{}
        T_k (\L) \le 2^{9k+4d} k^k |\L|^k \cdot 2^{ \frac{2s k (\log k)^2}{\log (k^{2s} |\L|^{s-2})}}
                        =  2^{9k+4d} k^k |\L|^k \cdot 2^{ \frac{6k (\log k)^2}{\log (k^{6} k^\_phi)}}
                            = 2^{9k+4d} k^k |\L|^k k^{\frac{6}{6+\_phi}}\,.
    \end{equation}
    On the other hand, using Theorem  \ref{t:main}, we get
    $T_k (\L) \ge \d \a^{2k} |\L|^{2k} / (2^{4k} \d^{2k})$.
    Whence $|\L| \le 2^{17 + 8d (\log (1/\d))^{-1}} (\d/\a)^2 \log^{\_phi} (1/\d)$ and
    $|\L^*| \le 2^{20 + 8d (\log (1/\d))^{-1}} (\d/\a)^2 \log^{\_phi} (1/\d)$.
    This completes the proof.



\refstepcounter{section} \label{finite_fields}

{\bf \arabic{section}. Some examples of sets of large exponential
sums in vectors spaces over finite field.}

Let $p$ be a prime number, $n$ and $N$ be positive integers,
$N=p^n$. In the section we consider groups $\Z_p^n = (\Z /
p\Z)^n$, $|\Z_p^n| = N$. A finite Abelian group $\Z_p^n$ is a
vector space with inner product
$$
    \v{x} \cdot \v{y} = <\v{x},\v{y}> = x_1 y_1 + \dots + x_n y_n \pmod p \,.
$$
Let $f: \Z_p^n \to \C$ be an arbitrary function. Denote by $\F{f}$
the Fourier transform of $f$
$$
    \F{f}(\v{r}) = \sum_{\v{x} \in \Z_p^n} f(\v{x}) e( - (\v{r}\cdot \v{x}) ) \,,
$$
where $e(x) = e^{2\pi i \frac{x}{p}}$, $x\in \Z_p$.

Let $\v{v}_1, \dots, \v{v}_k$ be linear--independent vectors and
let $\eps_1,\dots, \eps_k$ be elements of $\Z_p$. Define the {\it
affine subspace $P$} (of codimension $k$) by
$$
    P = P_{\eps_1,\dots, \eps_k} = \{ \v{x} \in \Z_p^n ~:~ <\v{x},\v{v}_1> = \eps_1, \dots, <\v{x},\v{v}_k> = \eps_k\} \,.
$$
It is easy to calculate the the Fourier transform of $P$. Let $L$
be the subspace of dimension $k$ spanned by
$\v{v}_1,\dots,\v{v}_k$. Suppose that $\v{r}\in \Z_p^n$ is an
arbitrary vector.
We have $\v{r} = \sum_{i=1}^k r_i \v{v}_i + \v{v}$, where $\v{v}
\in L^\bot$. Then
\begin{equation}\label{f:aff_subspace}
    \F{P} (\v{r}) = L(\v{r}) |P| \cdot e(-\sum_{i,j=1}^k \eps_i r_j <\v{v}_i, \v{v}_j> ) \,.
\end{equation}
Thus $|\F{P} (\v{r})|$ either equals zero or equals $|P|$.

In this section we consider the case $p=2$.
At the case the Fourier transform of a function $f$, $f: \Z_2^n
\to \C$ is equal to
$$
    \F{f}(\v{r}) = \sum_{\v{x} \in \Z_p^n} (-1)^{<\v{r}, \v{x}>} f(\v{x}) \,.
$$

First of all let us prove an analog of Theorem \ref{t:T2_low_Z_N}
for $\Z_2^n$.
It is very convenient to split our results into                                                         
Theorem \ref{t:T2_low'} and Theorem \ref{t:T2_low}. Theorem
\ref{t:T2_low'} is simpler then Theorem \ref{t:T2_low} but we need
in rigid condition (\ref{f:burden}) in our proof.

\Th \label{t:T2_low'} {\it
    Let $\d,\a \in (0,1]$ be real numbers, $\a \le \d/2$, $\d \le 2^{-5}$, and
    \begin{equation}\label{f:burden}
        \frac{2\d}{\a} \log \frac{1}{2\a} \le \log N \,.
    \end{equation}
    Then there exists a set $A\subseteq \Z_2^n$ such that  $\d N \le |A| \le 8\d N$,
    $|\r_\a (A)| \ge \frac{\d}{8\a^2}$ and for all $k$, $2 \le k \le 2^{-1} \log (1/8\d)$,
    we have
    $T_k (\r_\a (A)) \le \frac{8\d}{\a^{2k}}$.
}
\\
\Proof
    Let $\v{e}_1 = (1,0,\dots,0)$, $\v{e}_2 = (0,1,0,\dots,0), \dots,
    \v{e}_n = (0,\dots,0,1)$ be the standard basis of $\Z_2^n$.
    Let also $k' = [\log 1/(2\a)]$, $t=\lceil \d/\a \rceil$, $n=\log N$.
    Let $i\in [t]$ and let $P_i$ be a affine subspace such that
    $$
        P_i = \{ \v{x} \in \Z_2^n ~:~ <\v{x}, \v{e}_j> = 0\,, \quad j = (i-1)k' +1, \dots, ik' \} \,.
    $$
    Since $t k' \le \frac{2\d}{\a} \log \frac{1}{2\a} \le \log N = n$,
    it follows that all affine subspaces  $P_i$
    are well defined.
    Let $A = \bigcup_{i=1}^t P_i$.
    Clearly, $|A| \le t 2^{-k'} N \le 8\d N$.
    Let us prove that $|A| \ge \d N$.
    We have $|P_i| = N 2^{-k'}$, $i\in [t]$.
    Besides for any $l\in [t]$ and for all  {\it different} subspaces $P_{i_1}, \dots, P_{i_l}$
    the following holds
    \begin{equation}\label{f:aff_intersect}
        | P_{i_1} \cap \dots \cap P_{i_l} | = N 2^{-k'l} \,.
    \end{equation}
    Using the inequality  $\d \le 2^{-5}$, we get
    $$
        |A| \ge \sum_{i=1}^t |P_i| - \sum_{i,j=1, i\neq j}^t |P_i \cap P_j| \ge t 2^{-k'} N - t^2 (2^{-k'})^2 N
                = t 2^{-k'} N \left( 1 - \frac{t}{2^{k'}} \right) \ge \d N \,.
    $$
    Let us prove now that $|\r_\a (A)| \ge \frac{\d}{8\a^2}$.
    Let $L_i$ be a subspace of $\Z_2^n$ of the dimension $k'$ spanned by $\{ \v{e}_j \}_{,\,\, j = (i-1)k' +1, \dots, ik'}$.
    Suppose that $\v{s} \in \Z_2^n$ is an arbitrary vector.
    Using (\ref{f:aff_subspace}), we obtain
    \begin{equation}\label{f:TwO'}
        \F{P}_i (\v{s}) = |P_i| L_i (\v{s}) \,.
    \end{equation}
    Whence $\r_\a (A) \subseteq \bigcup_{i=1}^t L_i$.
    Prove that $\bigcup_{i=1}^t L_i \subseteq \r_\a (A)$.
    Obviously,  $\v{0} \in \r_\a (A)$.
    Let $\v{s}$ be a non--zero vector such that $\v{s}$ belongs to some $L_i$.
    Clearly, for any $i,j \in [t]$, $i\neq j$, we have $L_i \cap L_j = \{ \v{0} \}$.
    Using  this fact and (\ref{f:TwO'}), we get
    \begin{equation}\label{tmp:19:59_Aug18}
        \F{A} (\v{s})
            =
                \F{P}_i (\v{s})
                    -
                        \sum_{j=1}^t (P_i \cap P_j) \F{}\,\, (\v{s})
                            +
                                \sum_{j,l=1,\,\, j\neq l,\,\,\, j,l \neq i}^t (P_i \cap P_j \cap P_l) \F{}\,\, (\v{s})
                                    +
                                        \dots
    \end{equation}
    Using (\ref{f:aff_intersect}) and (\ref{tmp:19:59_Aug18}), we obtain
    \begin{equation}\label{}
        | \F{A} (\v{s}) |
            \ge
                2^{-k'} N - 2^{-k'} N \left( \frac{t}{2^{k'}} + \frac{t^2}{(2^{k'})^2} + \dots \right)
                    \ge
                        2^{-k'-1} N
                            \ge
                                \a N \,.
    \end{equation}
    Hence $\bigcup_{i=1}^t L_i \subseteq \r_\a (A)$
    and
    $|\r_\a (A)| \ge \sum_{i=1}^t |L_i| - t \ge t 2^{k'} - t \ge \frac{\d}{8\a^2}$.

    Finally, let us prove that for all $2\le k \le 2^{-1} \log (1/8\d)$,
    we have $T_k (\r_\a (A)) \le \frac{8\d}{\a^{2k}}$.
    Consider the equation
    \begin{equation}\label{f:TeN'}
        r_1 + \dots + r_k = r'_1 + \dots + r'_k \,,
    \end{equation}
    where all vectors $r_j$, $r'_j$ belong to $\r_\a (A)$.
    As was noted above $\r_\a (A) \subseteq \bigcup_{i=1}^t L_i$.
    Hence any vector in (\ref{f:TeN'}) belongs to some subspace $L_{i_j}$.
    Let $z$ be a non--negative integer,
    and
    $s_1, \dots, s_l$ be positive integers such that
    $s_1 + \dots + s_l + z = 2k$.
    By $E(s_1,\dots,s_l,z)$ denote the set of all solutions
    $r_1,\dots, r_k$, $r'_1,\dots, r'_k$ of (\ref{f:TeN'}) such that
    among $r_j$, $r'_j$
    there exist exactly $z$ of zeroes,
    there exist exactly $s_1$ non--zero residuals belong to a subspace $L_{j_1}$,
    there exist exactly $s_2$ non--zero residuals belong to a subspace $L_{j_2}$,
    $\dots$,
    there exist exactly $s_l$ non--zero residuals belong to a subspace $L_{j_l}$
    and at the same time all sets
    $L_{j_1}, L_{j_2}, \dots, L_{j_l}$ are different.
    We have
    $$
        T_k (\r_\a (A)) = \sum_{l=1}^{2k} \sum_{z=0}^{2k}~ \sum_{s_1,\dots,s_l,\,\, s_1+\dots+s_l+z = 2k}
                                |E(s_1,\dots,s_l,z)|
                        =
    $$
    \begin{equation}\label{f:-10}
                        =
                                t(2^{k'})^{2k-1}
                                    +
                                        \sum_{l=2}^{2k} \sum_{z=0}^{2k}~ \sum_{s_1,\dots,s_l,\,\, s_1+\dots+s_l+z = 2k}
                                            |E(s_1,\dots,s_l,z)| \,.
    \end{equation}
    Let us fixed $s_1, \dots, s_l,z$ and consider the solutions of (\ref{f:TeN'})
    belong to fixed subspaces $L_{j_1}, \dots, L_{j_l}$.
    Denote by $E(s_1,\dots,s_l,z) (L_{j_1}, \dots, L_{j_l})$
    the set of
    all these solutions.
    Rewrite (\ref{f:TeN'}) as
    \begin{equation}\label{}
        \v{u}_1 + \dots + \v{u}_l = \v{0} \,,
    \end{equation}
    where $\v{u}_i \in L_{j_i}$, $i\in [l]$.
    For all $i,h \in [t]$, $i\neq h$, we have
    $L_{j_i} \cap L_{j_h} = \{ \v{0} \}$.
    Hence all vectors $\v{u}_i$ equal $\v{0}$.
    It follows that
    $$
        | E(s_1,\dots,s_l,z) (L_{j_1}, \dots, L_{j_l}) |
            \le
                \frac{(2k)!}{s_1! \dots s_l! z!} (2^{k'})^{s_1-1} \m \dots \m (2^{k'})^{s_l-1}
                    \le
                         \frac{(2k)!}{s_1! \dots s_l! z!} (2^{k'})^{2k - l} \,.
    $$
    Whence
    \begin{equation}\label{f:-11}
        |E(s_1,\dots,s_l,z)|
            \le
                \binom{t}{l} \frac{(2k)!}{s_1! \dots s_l! z!} (2^{k'})^{2k - l}
                    \le
                        \frac{t^l}{l!} \cdot \frac{(2k)!}{s_1! \dots s_l! z!} (2^{k'})^{2k - l} \,.
    \end{equation}
    Combining (\ref{f:-11}) and (\ref{f:-10}), we get
    $$
        T_k (\r_\a (A))
            \le
                t(2^{k'})^{2k-1}
                    +
                        \sum_{l=2}^{2k} \frac{t^l}{l!} (2^{k'})^{2k - l}
                            \sum_{z=0}^{2k}~ \sum_{s_1,\dots,s_l,\,\, s_1+\dots+s_l+z = 2k} \frac{(2k)!}{s_1! \dots s_l! z!}
            \le
    $$
    \begin{equation}\label{}
            \le
                t(2^{k'})^{2k-1}
                    +
                        \sum_{l=2}^{2k} \frac{t^l}{l!} (2^{k'})^{2k - l} (l+1)^{2k}
            =
                t(2^{k'})^{2k-1}
                    +
                        (2^{k'})^{2k} \sum_{l=2}^{2k} \left( \frac{t}{2^{k'}} \right)^l \cdot (l+1)^{2k} \cdot \frac{1}{l!} \,.
    \end{equation}
    Consider the function $f(l) = (t/2^{k'})^l (l+1)^{2k}$.
    It is easy to see that $f(l)$ has maximum at $l_0 = 2k / \ln (2^{k'} / t) - 1$
    and for all $l\ge l_0$ the function $f(l)$ is monotonically decreasing.
    By assumption $k \le 2^{-1} \log (1/8\d)$.
    Hence $l_0 \le 1$.
    It follows that
    $$
        T_k (\r_\a (A))
            \le
                t(2^{k'})^{2k-1} + 2^{2k} t(2^{k'})^{2k-1}
                    \le
                        2^{2k+1} t(2^{k'})^{2k-1}
                            \le
                                2^{2k+1} \cdot \frac{2\d}{\a} \left( \frac{1}{2\a} \right)^{2k-1}
                                    =
                                        \frac{8\d}{\a^{2k}} \,.
    $$
    This completes the proof.

\Note In special cases of choosing $\d$ and $\a$ we do not need in
a bound $k \ll \log (1/\d)$ of Theorem \ref{t:T2_low'}. For
example, suppose that $\a \approx \d$ and $A$ is a subspace of
$\Z_2^n$ of codimension $k'$, $k' \approx \log ( 1/\d )$. Then
$\r_\a (A)$  is a subspace of dimension $k'$ and it has the
cardinality $\approx 1/ \d$. It is easy to see that for all $k\ge
2$, we have $T_k (\r_\a (A)) = ( 2^{k'} )^{2k-1} \approx 1 /
\d^{2k-1}$. This quantity coincides with lower bound
(\ref{f:T_k}).

We shall consider the simplest case of $k=2$ in our next Theorem
\ref{t:T2_low}, i.e. we shall prove that the lower bound for $T_2
(\r_\a (A))$ from Theorem \ref{t:main} is best possible.

\Th \label{t:T2_low} {\it
    Let $\d,\a \in (0,1]$ be real numbers, $32 \d^2 \le \a \le \d/2$, $\a \ge N^{-2^{-300}}$,
    $\a \le 2^{-30}$, and
    $\frac{\d}{\a} \log \frac{1}{2\a} \ge 400 \log N \cdot \log (8\log N)$.
    Then there exists a set $A\subseteq \Z_2^n$ such that  $\d N \le |A| \le 8\d N$,
    $|\r_\a (A)| \ge \frac{\d}{8\a^2}$ and $T_2 (\r_\a (A)) \le \frac{16\d}{\a^4}$.
}

To prove such result we need in a well--known large deviations
inequality of Bernstein \cite{Bernstein}. The following variant of
this inequality can be found in \cite{GreenA+A}.

\Th \label{t:large_dev} {\it
    Let $X_1,\dots,X_n$ be independent random variables with $\mathbb{E} X_j = 0$
    and $\mathbb{E} |X_j|^2 = \sigma_j^2$.
    Let $\sigma^2 = \sigma_1^2 + \dots + \sigma_n^2$.
    Suppose that for all  $j\in [n]$, we have $|X_j| \le 1$.
    Let also $t$ be a real number such that  $\sigma^2 \ge 6 nt$.
    Then
    $$
        \mathbb{P} \left( \left| \frac{X_1+\dots + X_n}{n} \right| \ge t \right) \le 4 e^{-n^2 t^2 / 8\sigma^2} \,.
    $$
}

Using Theorem \ref{t:large_dev}, we prove a combinatorial lemma.

\Lemma \label{l:family} {\it
    Let $n,k,r,t$ be real numbers, $4\le r \le k/2$, $2k\le n$,
    and
\begin{equation}\label{f:l_cond1+1}
    kt > 288 n \ln (8n) \quad \mbox{ and } \quad  t^2 \cdot \frac{2^k \binom{n-k}{k-\lceil k/r \rceil}}{\binom{n}{k}} \le 1/2 \,.
\end{equation}
    Then  there exist sets  $A_1, \dots, A_t \subseteq [n]$, $|A_i| = k$, $i\in [t]$
    such that \\
$1)~$ For all $i,j\in [t]$, $i\neq j$, we have $|A_i \cap A_j| < k/r$. \\
$2)~$ For any $i\in [t]$ there exist at most $2tk^2/n$
      sets $A_j$ such that $A_j \cap A_i \neq \emptyset$.
}
\\
\Proof
    Let $\Omega$ be a family of all subsets of $[n]$ of the cardinality $k$,
    $|\Omega| = \binom{n}{k} = M$.
    Choose sets $A_1, \dots, A_t \in \Omega$ at random
    (uniformly and independently).

    Let $U_{ij}$, $i,j \in [t]$, $i\neq j$ be a random event consists in
    $|A_i \cap A_j| \ge k/r$.
    Let also $U = \bigcup_{i,j\in [t], i\neq j} U_{ij}$.
    Let us fix a set $A_i$.
    It is easy to see that there are exactly
    $$
        \sigma := \sum_{l=\lceil k/r \rceil}^k \binom{k}{l} \binom{n-k}{k-l}
    $$
    sets $A_j \in \Omega$ such that $|A_j \cap A_i| \ge k/r$.
    Hence the probability of $U_{ij}$ is equal to $\sigma / M$.
    Whence $\mathbb{P} (U) \le t^2 \sigma / M$.

    Let $x \in [n]$ and
    $\xi_j^x (\o)$, $\o \in \Omega^t$, $j\in [t]$ be a random variable such that
    $\xi_j^x (\o) = 1$ if $x\in \o_j$ and  $\xi_j^x (\o) = 0$
    otherwise.
    Clearly, $\xi^x (\o) := \sum_{j=1}^t \xi_j^x (\o)$ is the number of sets $A_j$
    such that $x \in A_j$.
    Further, for any $x$ and $j$, we have $\mathbb{E} \xi_j^x = k/n$ and $\mathbb{D} \xi_j^x = k/n - (k/n)^2$.
    Let $x\in [n]$ and let $V_x$ be the event such that $x$ belongs to at least
    $7tk/(6n)$ the sets $A_j$.
    Let also $V = \bigcup_{x\in [n]} V_x$.
    Using Theorem \ref{t:large_dev}, we get
    $$
        \mathbb{P} (V_x) \le
                            \mathbb{P} \left( \o ~:~ \left| \xi^x (\o) - \frac{tk}{n} \right| > \frac{tk}{6n} \right)
                                \le 4 e^{-kt/(288n)} \,.
    $$
    Whence
    \begin{equation}\label{f:tmp:20:05_Aug_15}
        \mathbb{P} (V) \le \sum_{x\in [n]} \mathbb{P} (V_x) \le 4 n e^{-kt/(288n)} \,.
    \end{equation}
    By assumption $kt > 288 n \ln (8n)$.
    It follows that $4 n e^{-kt/(288n)} < 1/2$.
    Besides $\sigma \le 2^k \binom{n-k}{k-\lceil k/r \rceil}$.
    Using condition (\ref{f:l_cond1+1}) and inequality (\ref{f:tmp:20:05_Aug_15}), we obtain
    $$
        \mathbb{P} (U\cup V) \le \mathbb{P} (U) + \mathbb{P} (V) \le t^2 \frac{\sigma}{M} + 4 n e^{-kt/(288n)}
        \le t^2 \cdot \frac{2^k \binom{n-k}{k-\lceil k/r \rceil}}{\binom{n}{k}} + 4 n e^{-kt/(288n)}
        < 1/2 + 1/2 = 1\,.
    $$
    Hence there exists a collection of sets $A_1,\dots,A_t \subseteq [n]$, $|A_j| = k$
    such that $1)$ holds and such that any element $x\in [n]$ belongs to at most
    $7tk/(6n) \le 2tk/n$ of these sets.
    Using the last inequality, we get that for all $i\in [t]$ there are
    at most $2tk^2/n$ the sets $A_j$ such that $A_i\cap A_j \neq \emptyset$.
    This concludes the proof.

{\bf Proof of Theorem \ref{t:T2_low}.}
    Let $\v{e}_1 = (1,0,\dots,0)$, $\v{e}_2 = (0,1,0,\dots,0), \dots,
    \v{e}_n = (0,\dots,0,1)$ be the standard basis of $\Z_2^n$.
    Let also $r=32$, $k = [\log 1/(2\a)]$, $t=\lceil \d/\a \rceil$, $n=\log N$.
    By assumption $\frac{\d}{\a} \log \frac{1}{2\a} \ge 400 \log N \cdot \log (8\log N)$.
    Hence $kt > 288 n \ln (8n)$.
    Besides $\a \ge N^{-2^{-300}}$.
    Whence
    $$
        t^2 \cdot \frac{2^k \binom{n-k}{k - \lceil k/32\rceil}}{\binom{n}{k}} \le t^2 2^k \frac{k^{k/32 + 1}}{(n-k)^{k/32}}
        \le k t^2 2^k \left( \frac{2k}{n} \right)^{k/32} \le 1/2 \,.
    $$
    Using Lemma \ref{l:family}, we find a collection of sets $A_1,\dots, A_t$
    such that $1)$ and $2)$ hold.

    We shall construct a family of affine subspaces $P_1, \dots, P_t$
    of such form
    $$
        P_i = P_i^{\v{\eps}} = \{ \v{x} \in \Z_2^n ~:~ <\v{x}, \v{e}_j> = \eps_i^{(j)}\,, \quad j\in A_i \} \,,
    $$
    where $\v{\eps}_i = ( \eps_i^{(j)})$ be a vector from $\Z_2^k$.
    Thus to construct affine subspaces $P_i$, we need to
    choose vectors $\v{\eps}_1,\dots, \v{\eps}_t$.
    Let $\v{\eps}_1 = \v{0}$ and we obtain $P_1$.
    Suppose that we have the affine subspaces $P_1,\dots,P_d$.
    Let us construct a vector $\v{\eps}_{d+1}$ and a affine subspace $P_{d+1}$.
    Let $C_d = \bigcup_{i=1}^d P_i$.
    Clearly, $|C_d| \le d N 2^{-k} \le t N 2^{-k} \le 8\d N$.
    Let $\v{\eps}_{d+1}$ be a vector such that
    \begin{equation}\label{tmp:23:09_16Aug}
        |P_{d+1}^{\v{\eps}_{d+1}} \bigcap C_{d}| \le 2\d \cdot 2^{-k} N \,.
    \end{equation}
    Since
    $$
        \sum_{\v{\eps} = ( \eps_i^{(j)})\,,~ j\in A_r}  |C_d \cap P_{d+1}^{\v{\eps}}| = |C_d| \,,
    $$
    it follows that such a vector $\v{\eps}_{d+1}$ exists.
    So we have the affine subspaces $P_1,\dots,P_t$.
    Let $A = C_t = \bigcup_{i=1}^t P_i$.
    Clearly, $|A| \le 8\d N$.
    Let us prove that $|A| \ge \d N$.
    We have $|P_i| = N 2^{-k}$, $i\in [t]$.
    Using (\ref{tmp:23:09_16Aug}), we get
    \begin{equation}\label{f:expl_1}
        |A| = |C_t| = |C_{t-1}| + |P_t| - |C_{t-1} \cap P_t| \ge |C_{t-1}| + N 2^{-k} - \frac{8\d N}{2^k}
        \ge
    \end{equation}
    \begin{equation}\label{f:expl_2}
        \ge
            |C_{t-2}| + 2N 2^{-k} - 2 \frac{8\d N}{2^k} \ge \dots \ge t N 2^{-k} - t \frac{8\d N}{2^k}
            = tN 2^{-k} (1-8\d) \ge \d N \,.
    \end{equation}

    Let us prove that $|\r_\a (A)| \ge \frac{\d}{8\a^2}$.
    Let $L_i$ be a subspace of $\Z_2^n$ of the dimension $k$ spanned by $\{ \v{e}_j \}_{j\in A_i}$.
    Let also
    $$
        M_i  = \{ \v{x} \in L_i ~:~ \mbox{ the number of units in } \v{x} \mbox{ at least } k/8 \} \,.
    $$
    Clearly, for all $i\in [t]$, we get $|M_i| \ge 2^{k-1}$.
    Since for any $i,j \in [t]$, $i\neq j$, we have  $|A_i \cap A_j| < k/r < k/8$,
    it follows that for all $i,j \in [t]$, $i\neq j$, we obtain $M_i \cap L_j = \emptyset$.
    In particular, for any $i,j \in [t]$, $i\neq j$, we get $M_i \cap M_j = \emptyset$.
    Let $\v{s} \in \Z_2^n$ be an arbitrary vector.
    Using (\ref{f:aff_subspace}), we get
    \begin{equation}\label{f:TwO}
        \F{P}_i (\v{s}) = e(-\sum_{j\in A_i} \eps_i^{(j)} s_j) |P_i| L_i (r) \,.
    \end{equation}
    Whence $\r_\a (A) \subseteq \bigcup_{i=1}^t L_i$.
    Prove that $\bigsqcup_{i=1}^t M_i \subseteq \r_\a (A)$.
    Let $i\in [t]$, and $\v{s} \in M_i$ be a vector.
    We have
    $$
        \F{A} (\v{s}) = \F{P}_t (\v{s}) + \F{C}_{t-1} (\v{s}) + \theta \frac{8\d N}{2^k} \,,
    $$
    where $|\theta| \le 1$.
    By the same arguments as in (\ref{f:expl_1}) --- (\ref{f:expl_2}), we get
    \begin{equation}\label{f:ToW_primes}
        \F{A} (\v{s}) = \sum_{l=1}^t \F{P}_l (\v{s}) + \t{\theta} t \frac{8\d N}{2^k} \,,
    \end{equation}
    where $|\t{\theta}| \le 1$.
    For all $i,j \in [t]$, $i\neq j$, we have $M_i \cap L_j = \emptyset$.
    Hence
    \begin{equation}\label{f:ThRee}
        | \sum_{l=1}^t \F{P}_l (\v{s}) | = | \F{P}_i (\v{s}) | = N 2^{-k} \,.
    \end{equation}
    Using (\ref{f:ToW_primes}), (\ref{f:ThRee}) and $\a \ge 32 \d^2$, we obtain
    $$
        | \F{A} (\v{s}) | \ge N 2^{-k} - t \frac{8\d N}{2^k} = \frac{N}{2^k} ( 1 - 8\d t ) \ge \a N \,.
    $$
    Therefore $\bigsqcup_{i=1}^t M_i \subseteq \r_\a (A)$ and
    $|\r_\a (A)| \ge t 2^{k-1} \ge \frac{\d}{8\a^2}$.

    Finally, we shall show that $T_2 (\r_\a (A)) \le \frac{16\d}{\a^4}$.
    Consider the equation
    \begin{equation}\label{f:TeN}
        \v{r}_1 + \v{r}_2 = \v{r}_3 + \v{r}_4 \,,
    \end{equation}
    where all $\v{r}_l$, $l = 1,2,3,4$ belong to $\r_\a (A)$.
    As was noted above $\r_\a (A) \subseteq \bigcup_{i=1}^t L_i$.
    It follows that any vector $\v{r}_l$ belongs to some
    subspace $L_{i_l}$.
    Let $M = \bigsqcup_{i=1}^t M_i$, and $Q = (\bigcup_{i=1}^t ) \setminus M$.
    For any $i\in [t]$, we have
    \begin{equation}\label{}
        |L_i \setminus M_i| = \sum_{l=1}^{[k/8]} \binom{k}{l} \le \frac{k}{8} \binom{k}{[k/8]} \,.
    \end{equation}
    Hence
    \begin{equation}\label{f:-2}
        |Q| \le \sum_{i=1}^t |L_i \setminus M_i| \le \frac{kt}{8} \binom{k}{[k/8]} \,.
    \end{equation}
    Using Stirling's formula, (\ref{f:-2}) and $\a \ge 8 \d^2$, we obtain
    \begin{equation}\label{}
        T_2 (Q) \le |Q|^3 \le \frac{k^3 t^3}{8^3}  \binom{k}{[k/8]}^3 \le \frac{t}{8} (2^k)^3 \,.
    \end{equation}
    The last inequality implies that
    $$
        T_2 (\r_\a (A)) \le T_2 ( M\sqcup Q) = \frac{1}{N} \sum_{\v{r} \in \Z_2^n} |\F{M} (r) + \F{Q} (r)|^4
            \le \frac{8}{N} \sum_{\v{r} \in \Z_2^n} |\F{M} (r)|^4
                +
                \frac{8}{N} \sum_{\v{r} \in \Z_2^n} |\F{Q} (r)|^4
                    \le
    $$
    \begin{equation}\label{f:-4}
                    \le 8 T_2(M) + 8 T_2 (Q) \le 8 T_2 (M) + t (2^k)^3 \,.
    \end{equation}
    Thus to obtain an upper bound for
    $T_2 (\r_\a (A))$
    we need to compute $T_2 (M)$.

    So let $\v{r}_l \in M_{i_l}$, $l\in [4]$.
    For all $i,j \in [t]$, $i\neq j$, we have $|A_i \cap A_j| < k/r$.
    Since $3k/r = 3k/32 < k/8$, it follows that all set $M_{i_l}$, $l=1,2,3,4$
    cannot be different.
    Furthermore we have three cases : \\
    $1)~$ $i_1 = i_2 = i_3 = i_4\,,$ \\
    $2)~$ $i_1 = i_3$,  $i_2 = i_4$ and $i_1\neq i_2\,,$ \\
    $3)~$ $i_1 = i_4$,  $i_2 = i_3$ and $i_1\neq i_2$. \\
    In the first case the number of solutions of (\ref{f:TeN}) does not exceed $t (2^k)^3$.
    Let us consider the case $2)$ (or $3)$).
    Let us fixed $i_1$ and $i_2$, $i_1 \neq i_2$.
    Let $\v{u} = \v{r}_1 - \v{r}_3 = \v{r}_4 - \v{r}_2$.
    Clearly, $\v{u} \in L_{i_1} \cap L_{i_2}$.
    If $A_{i_1} \cap A_{i_2} = \emptyset$ then $\v{u} = \v{0}$ and
    $\v{r}_1 = \v{r}_3$, $\v{r}_2 = \v{r}_4$.
    Hence if $A_{i_1} \cap A_{i_2} = \emptyset$ then (\ref{f:TeN})
    has at most $(2^k)^2$ solutions.
    Suppose that $A_{i_1} \cap A_{i_2} \neq \emptyset$.
    Since $|A_{i_1} \cap A_{i_2}| < k/r$, it follows that
    the number of solutions of (\ref{f:TeN})
    does not exceed $(2^k)^2 \cdot 2^{k/r}$ in this case.
    Whence the number of solutions of (\ref{f:TeN}) at most
    $$
        \sum_{i_1=1}^t ~ \sum_{i_2=1, i_2 \neq i_1, A_{i_1} \cap A_{i_2} = \emptyset}^t (2^k)^2
            +
        \sum_{i_1=1}^t ~ \sum_{i_2=1, i_2 \neq i_1, A_{i_1} \cap A_{i_2} \neq \emptyset}^t (2^k)^2 \cdot 2^{k/r}
            := \sigma_1 \,.
    $$
    Using property $2)$ of the collection of the sets $A_1,\dots, A_t$, we obtain
    that the number of $i_2 \neq i_1$
    such that $A_{i_1} \cap A_{i_2} \neq \emptyset$ does not exceed $2tk^2/n$.
    Hence
    $$
        \sigma_1 \le t^2 (2^k)^2 + t \frac{2tk^2}{n} (2^k)^2 \cdot 2^{k/32} \,.
    $$
    Using the last inequality and $\a \ge 32 \d^2$, we get
    $\sigma_1 \le t(2^k)^3 + t (2^k)^3 = 2t (2^k)^3$.
    Whence the {\it total} number of solutions of (\ref{f:TeN}) at most
    $$
        T_2 (\r_\a (A)) \le 8 (t (2^k)^3 + 2t (2^k)^3 + 2t (2^k)^3) + t (2^k)^3 = 41 t (2^k)^3
            \le 41 \frac{2\d}{\a} \frac{1}{(2\a)^3} \le \frac{16\d}{\a^4} \,.
    $$
    This completes the proof.


    As was showed in Theorems \ref{t:T2_low'}, \ref{t:T2_low} an upper bound
    of Theorem \ref{t:main} is best possible.
    It is easy to see that the number of elements $\l^*_i$ in (\ref{f:r=log})
    of Theorem \ref{t:Chang_log}
    is also best possible.
    Indeed, let $\a \approx \d$ and let $A$ be a subset of $\Z_2^n$
    such that $|\r_\a (A)| \approx \d /\a^2 \approx 1/\d$.
    Certainly, such sets exist, for example one can take a subspace of
    $\Z_2^n$ of the cardinality $\d N$.
    By Chang's Theorem
    there exists a set $\L^*$, $|\L^*| \ll \log (1/\d)$ such that
    for any $\v{r} \in \r_\a (A)$, we have (\ref{f:r=log}).
    On the other hand since $|\r_\a (A)| \approx 1/\d$, it follows that
    there exists a vector $\v{r} \in \r_\a (A)$ such that we need in $k \gg \log (1/\d)$ vectors of $\L^*$
    to involve $\v{r}$ in some equation (\ref{f:r=log}).
    Indeed, we have at most $2^k \binom{|\L^*|}{k}$ of linear combinations
    of $k$ vectors from $\L^*$.
    Hence the following inequality must be hold $2^k \binom{|\L^*|}{k} \gg 1/\d$.
    This implies that $k \gg \log (1/\d)$.


\begin{thebibliography}{99}




    \bibitem{Gow_surv} \emph{Gowers\,W.\,T. }
    Rough structure and classification //
    Geom. Funct. Anal., Special Volume - GAFA2000 "Visions in Mathematics", Tel
    Aviv, (1999) Part I, 79--117.



    \bibitem{Gow_m} \emph{Gowers\,W.\,T. }
    A new proof of Szemer\'{e}di's theorem //
    Geom. Funct. Anal. \textbf{11} (2001), 465--588.



    \bibitem{Ch_Fr} \emph{Chang\,M.--\,C., }
    A polynomial bound in Freiman's theorem //
    Duke Math. J. \textbf{113} (2002) no. 3, 399--419.



    \bibitem{Ruzsa_Freiman} \emph{Ruzsa\,I. }
    Generalized arithmetic progressions and sumsets //
    Acta Math. Hungar., \textbf{65} (1994), 379--388.



    \bibitem{Bilu} \emph{Bilu\,Y. }
    Structure of sets with small sumset //
    Structure Theory of Sets Addition, Ast\'{e}risque, Soc. Math. France, Montrouge, \textbf{258} (1999), 77--108.



    \bibitem{Freiman} \emph{Freiman\,G.\,A. }
    Foundations of a Structural Theory of Set Addition /
    Kazanskii Gos. Ped. Inst., Kazan, 1966.
    Translations of Mathimatical Monographs \textbf{37}, AMS, Providence, R.I., USA.




    \bibitem{GreenA+A} \emph{Green\,B. }
    Arithmetic Progressions in Sumsets //
    Geom. Funct. Anal., \textbf{12} (2002) no. 3, 584--597.



    \bibitem{Green_Chang_exact} \emph{Green\,B. }
    Some constructions in the inverse spectral theory of cyclic groups //
    Comb. Prob. Comp. \textbf{12} (2003) no. 2, 127--138.



    \bibitem{Green_Chang1} \emph{Green\,B. }
    Spectral structure of sets of integers //
    Fourier analysis and convexity (survey article, Milan 2001),
    Appl. Numer. Harmon. Anal., Birkhauser Boston, Boston, MA (2004), 83--96.



    \bibitem{Green_Chang2} \emph{Green\,B. }
    Structure Theory of Set Addition //
    ICMS Instructional Conference in Combinatorial Aspects of
    Mathematical Analysis, Edinburgh March 25 --- April 5 2002.


    \bibitem{Green_Sz_Ab} \emph{Green\,B. }
    A Szemer\'{e}di--type regularity lemma in abelian groups //
    Geom. Funct. Anal. \textbf{15} (2005) no. 2, 340--376.


    \bibitem{Green_finite_fields} \emph{Green\,B. }
    Finite field model in additive combinatorics //
    Surveys in Combinatorics 2005, LMS Lecture Notes \textbf{329}, 1--29.



    \bibitem{Schoen} \emph{Schoen\,T. }
    Linear equations in $\Z_p$ //                                                 
    London Math. Soc., submitted.










    \bibitem{RuzsaA+A} \emph{Ruzsa\,I. }
    Arithmetic progressions in sumsets //
    Acta Arith. \textbf{60} (1991) no. 2, 191--202.


    \bibitem{Freiman_Yudin} \emph{Yudin\,A.\,A. }
    On the measure of large values of a trigonametric sum //
    Number Theory (under the edition of G.A. Freiman, A.M. Rubinov, E.V. Novosyolov),
    Kalinin State Univ., Moscow (1973), 163--174.



    \bibitem{Besser} \emph{Besser\,A. }
    Sets of integers with large trigonometric sums //
    Ast\'{e}risque \textbf{258} (1999), 35--76.


    \bibitem{Lev} \emph{Lev\,V.\,F. }
    Linear Equations over ${\mathbb F}_p$ and Moments of Exponential Sums //
    Duke Mathematical Journal \textbf{107} (2001), 239--263.


    \bibitem{Konyagin_Lev} \emph{Konyagin\,S.\,V., Lev\,V.\,F. }
    On the distribution of exponential sums //
    Integers: Electronic Journal of Combinatorial Number Theory \textbf{0} \# A01, (2000).


    \bibitem{Katznelson} \emph{de Leeuw\,K., Katznelson\,Y., Kahane\,J.\,P. }
    Sur les coefficients de Fourier des fonctions continues //
    C. R. Acad. Sci. Paris S\'{e}r. A--B \textbf{285} (1977) no. 16, A1001--A1003.


    \bibitem{Nazarov} \emph{Nazarov\,F.\,L. }
    The Bang solution of coefficient problem //
    Algebra i Analiz \textbf{9} (1997) no. 2, 272--287.
    English Transl. in St. Petersburg Math. J. \textbf{9} (1998) no. 2, 407--419.




    \bibitem{Ball} \emph{Ball\,K. }
    Convex geometry and functional analysis //
    Handbook of the geometry of Banach spaces, vol. I, North--Holland, Amsterdam (2001), 161--194.


    \bibitem{Rudin_book} \emph{Rudin\,W. }
    Fourier analysis on groups /
    Wiley 1990
    (reprint of the 1962 original).


    \bibitem{Rudin} \emph{Rudin\,W. }
    Trigonometric series with gaps //
    J. Math. Mech. \textbf{9} (1960), 203--227.


    \bibitem{Spencer_Beck} \emph{Spencer\,J. }
    Six Standard Deviations Suffice //
    Transactions of the American Mathematical Society \textbf{289} (1985), 679--706.


    \bibitem{Bernstein} \emph{Bernstein\,S. }
    Sur une modification de l'in\'{e}qualit\'{e} de Tchebichef //
    Annal. Sci. Inst. Sav. Ukr. Sect. Math. I (1924).


    \bibitem{Vinorgadov} \emph{Vinogradov\,I.\,M. }
    The method of trigonometric sums in number theory /
    M.: Nauka, 1971.


    \bibitem{Linnik} \emph{Linnik\,Y.\,V. }
    On Weyl's sums. //
    Math. Sbornik \textbf{12} (1943) I, 28--39.


    \bibitem{Nesterenko} \emph{Nesterenko\,Y.\,V. }
    On I.M. Vinogradov's mean--value theorem //
    Trudi of Moscow Math. Soc. \textbf{48} (1985), 97--105.


    \bibitem{Ruzsa_independent} \emph{Bajnok\,B., Ruzsa\,I. }
    The independence number of a subset of an abelian group //
    Integers: Electronic Journal of Combinatorial Number Theory \textbf{3} \# A02, 2003.                 


    \bibitem{Bu} \emph{Bourgain\,J. }
    On triples in arithmetic progression //
    Geom. Funct. Anal. \textbf{9} (1999), 968--984.


    \bibitem{Sh_dokl_exp1} \emph{Shkredov\,I.\,D. }
    On sets of large exponential sums //
    Doklady of Russian Academy of Sciences, 411, N 4,
    2006.


    \bibitem{Sh_exp1} \emph{Shkredov\,I.\,D. }
    On sets of large exponential sums //
    Izvestiya of Russian Academy of Sciences, submitted.










\end{thebibliography}
\end{document}